\documentclass[12pt]{article}

\usepackage{graphicx}
\usepackage{latexsym,amsmath,amsfonts,amscd, amsthm}
\usepackage{color}
\usepackage{changebar}
\usepackage{scalerel}
\graphicspath{{./}{./figure/}}

\newtheoremstyle{plainNoItalics}{}{}{\normalfont}{}{\bfseries}{.}{ }{}

\theoremstyle{plain}
\newtheorem{thm}{Theorem}[section]

\theoremstyle{plainNoItalics}

\newtheorem{defn}[thm]{Definition}
\newtheorem{rem}[thm]{Remark}
\newtheorem{prop}[thm]{Proposition}

\newcommand{\beq}{\begin{equation}}
\newcommand{\eeq}{\end{equation}}
\newcommand{\beqa}{\begin{eqnarray}}
\newcommand{\eeqa}{\end{eqnarray}}
\newcommand{\bit}{\begin{itemize}}
\newcommand{\eit}{\end{itemize}}
\newcommand{\bedef}{\begin{defn}}
\newcommand{\edefn}{\end{defn}}
\newcommand{\bpro}{\begin{prop}}
\newcommand{\epro}{\end{prop}}
\newcommand{\kv}{^{[k]}}

\newcommand{\kvm}{^{[k-1]}}

\newcommand{\kfm}{^{(k-1)}}

\def\Box{\mbox{ }\rule[0pt]{1.5ex}{1.5ex}}
\makeatletter

\newcommand{\Rmnum}[1]{\expandafter\@slowromancap\romannumeral #1@}
\makeatother

\begin{document}

\baselineskip=1.8pc


\begin{center}
{\bf
A High Order Time Splitting Method Based on Integral Deferred Correction
for Semi-Lagrangian Vlasov Simulations
}
\end{center}

\vspace{.2in}
\centerline{
Andrew Christlieb
\footnote{
Department of Mathematics, Michigan State University, East Lansing, MI, 48824. E-mail: andrewch@math.msu.edu
},
Wei Guo \footnote{Department of
Mathematics, University of Houston, Houston, TX, 77204. E-mail:
wguo126@math.uh.edu},
Maureen Morton
\footnote{ Department of Mathematics, Stark State College, North Canton, OH, 44270. E-mail: mmorton@starkstate.edu
},
Jing-Mei Qiu \footnote{Department of Mathematics, University of Houston,
Houston, TX, 77204. E-mail: jingqiu@math.uh.edu.
The second and the fourth authors are supported in part by Air Force Office of Scientific Computing YIP grant FA9550-12-0318, NSF grant DMS-1217008 and University of Houston.}
}

\bigskip
\noindent
{\bf Abstract.}
Semi-Lagrangian schemes with various splitting methods, and with different reconstruction/interpolation strategies have been applied to kinetic simulations. For example, the order of spatial accuracy of the algorithms proposed in {[Qiu and Christlieb, J.~Comp.~Phys., 2010]} is very high (as high as ninth order). However, the temporal error is dominated by the operator splitting error, which is second order for Strang splitting. It is therefore important to overcome such low order splitting error, in order to have numerical algorithms that achieve higher orders of accuracy in both space and time. In this paper, we propose to use the integral deferred correction (IDC) method to reduce the splitting error. Specifically, the temporal order accuracy is increased by $r$ with each correction loop in the IDC framework, where $r=1,\,2$ for coupling the first order splitting and the Strang splitting, respectively. The proposed algorithm is applied to the Vlasov-Poisson system, the guiding center model, and two dimensional incompressible flow simulations in the vorticity stream-function formulation.
We show numerically that the IDC procedure can automatically increase the order of accuracy in time. We also investigate numerical stability of the proposed algorithm via performing Fourier analysis to a linear model problem.
\vfill


\newpage

\newpage

\section{Introduction}
\label{sec1}
\setcounter{equation}{0}
\setcounter{figure}{0}
\setcounter{table}{0}

In this paper, we propose a semi-Lagrangian (SL) algorithm for Vlasov simulations, that is designed to be high order accurate in both space and time. In the proposed scheme, we utilize the high order weighted essentially non-oscillatory (WENO) spatial reconstructions \cite{Qiu_Christlieb, qiu_shu_sl, Qiu_Shu2}, coupled to a dimensional splitting SL framework, originally proposed in \cite{cheng1975integration}, with an integral deferred correction (IDC) framework \cite{DGR, highRKsdc}.

A simple collisionless model in describing plasma is given by the Vlasov-Poisson (VP) system,
 \begin{align}
 \label{eqn:vlasov3d}
	f_t
	+ {\bf v} \cdot \nabla_{\bf x} f 
	+ {\bf E}(t,{\bf x}) \cdot \nabla_{\bf v} f 
	&= 0, \\
	\label{eq: poisson}
	{\bf E}(t, {\bf x}) = - \nabla_{\bf x} \phi, \quad - \Delta_{\bf x} \phi &= -1+\rho(t,{\bf x}),
 \end{align}
 where $f (t, {\bf x}, {\bf v})$ describes the probability of finding a particle with velocity ${\bf v} \in \mathcal{R}^3$ at position ${\bf x}\in \mathcal{R}^3$ at time $t$, ${\bf E}$ is the electrostatic field, $\phi$ is the self-consistent electrostatic potential, and $\rho (t, {\bf x}) = \int{ f (t, {\bf x}, {\bf v})d{\bf v}}$ is the electron charge density and the $1$ represents the uniformly distributed infinitely massive ions in the background.  All physical constants in \eqref{eqn:vlasov3d} have been normalized to one.

Popular numerical approaches in solving the VP system can be classified as three types: Eulerian, Lagrangian and SL. The Lagrangian type particle methods, e.g. particle-in-cell (PIC) \cite{birdsall2005plasma, hockney2010computer, verboncoeur2005particle}, evolve the solution by following the trajectories of some sampled macro-particles, while the Eulerian approach  \cite{zhou2001numerical, cheng2012study, cheng2013energy} evolves the state variable according to the PDE on a fixed numerical grid. The SL approach is a mixed approach of Lagrangian and Eulerian in the sense that it has a fixed numerical grid; however, over each time step the state variable is evolved by propagating information along characteristics. The Eulerian and the SL approaches can be designed to be of very high order accuracy, an advantage when compared with the Lagrangian approach. On the other hand, because of the evolution mechanism, the SL method does not suffer the CFL time step restriction as in the explicit Eulerian approach, allowing for extra large time step evolution, and therefore less computational effort.

The SL approach is very popular in numerical weather prediction \cite{staniforth1991semi, Guo2013discontinuous} and kinetic simulations \cite{begue1999two, besse2003semi, huot2003instability, FilbetS, umeda2006comparison}, among many others. In kinetic simulations of plasma, a very popular approach is the Strang splitting SL method, first proposed by Cheng and Knorr \cite{cheng1975integration}. The advantage of performing such a splitting is that the decoupled lower dimensional equations of spatial advection and velocity ac/decelaration respectively are linear and are much easier to evolve numerically. Because of this, numerical schemes with high order of spatial accuracy have been designed, and demonstrated numerically to be more efficient than lower order schemes in finite difference framework with different reconstruction/interpolation strategies, such as, the cubic spline interpolation \cite{sonnendruecker}, the cubic interpolated propagation \cite{nakamura1999cubic}, the WENO interpolation \cite{carrillo2007nonoscillatory, Qiu_Christlieb, qiu_shu_sl}; in finite volume framework for the VP system \cite{filbet2001conservative} and for the guiding center Vlasov equation \cite{crouseilles2010conservative}; and in finite element discontinuous Galerkin framework \cite{James, SLDG}. On the other hand, the numerical error in time is dominated by the splitting error, which is relatively low order ($\mathcal{O}(\Delta t^2)$). In this paper, we will adopt the idea of integral deferred correction (IDC) to correct lower order dimensional splitting error.

The IDC methods are considered one-step, multi-stage integrators for solving initial value problems (IVPs). They are motivated by the defect/deferred correction (DC) methods \cite{stetter, FrUbTheoretical, pereyra}, and more recently the spectral deferred correction (SDC) methods \cite{DGR}. By construction, the IDC framework can systematically extend simple low order time integrators to high order ones by correcting provisional solutions.
The DC/SDC/IDC methods and their variants have been applied to many application areas such as chemical rate equations, reactive flow \cite{layton2004conservative, bourlioux2003high}, hyperbolic equations \cite{xiaxushu}, and parabolic equations \cite{christlieb2013high}. Additionally, recent developments in IDC algorithms have opened up new possibilities for increased computational speed via parallelization \cite{ridc2009}.

In this paper, we investigate applying the IDC framework to correct the dimensional splitting error in solving the VP system, and guiding center models with plasma applications, as well as in simulating incompressible flows. We choose the dimensional splitting SL finite difference WENO scheme \cite{Qiu_Christlieb, qiu_shu_sl} as the base solver in the prediction and correction steps of IDC, but point out that SL finite difference schemes with different reconstruction procedures such as those in \cite{sonnendruecker, nakamura2001exactly, carrillo2007nonoscillatory} can also be used. In the IDC framework, the low order temporal accuracy from dimensional splitting is increased by iteratively approximating error functions via solving error equations. In particular, the temporal error accuracy is lifted by $r$ in each correction loop, where $r=1,\,2$ for coupling the first order splitting and the Strang splitting, respectively. Our proposed SL scheme coupled with the IDC method enjoys the simplicity of the dimensional splitting algorithm, maintains the high order spatial accuracy, and can be designed to be of high order in time. However, the IDC methods with correction steps will render some CFL time step restrictions. We perform linear stability analysis via the classic Fourier approach and provide upper bounds of the CFL numbers for the proposed schemes. The CFL time step restriction is comparable to that for Eulerian methods using Runge-Kutta time discretization, leading to the computational cost at a similar scale.

We would like to remark that prior to the proposed IDC method in correcting the dimensional splitting error, several constructions of high order splitting methods have been developed. The methods proposed in \cite{yoshida, hairerGeom} are in the spirit of composition methods, requiring backward steppings (negative time steps). The number of intermediate stages, hence the computational cost, increases exponentially with the order of the splitting method. A fourth order splitting in \cite{yoshida} is applied to the VP system in \cite{James}. Higher order splitting methods that do not require backward steppings by using complex coefficients are proposed in \cite{highsplitanalySG}. The number of intermediate stages scales similarly to \cite{yoshida}. A fourth order splitting method for a linear Vlasov equation was presented in \cite{schaeffer}. However, the generalization of this fourth order method to a nonlinear problem is not straightforward.

The paper is organized as follows. Section~\ref{sec2} is a review of the dimensional splitting SL WENO method for the VP system. Section~\ref{sec3} first reviews the IDC method for IVPs in Section~\ref{sec3.1}, then formulates the error equations for the dimensional splitting algorithm for the VP system in Section~\ref{sec: 3.2}. The algorithm is then extended to guiding center models in Section~\ref{sec: 3.3}. Section~\ref{sec4} investigates the numerical stability via symbolic Fourier analysis for a linear model problem. Section~\ref{sec5} demonstrates the performance of the proposed schemes in simulating the VP system, the guiding center model, and advection in incompressible flow.
Finally, a conclusion is given in Section~\ref{conclusion}.

\section{Dimensional splitting SL WENO method for the VP system}
\label{sec2}
\setcounter{equation}{0}
\setcounter{figure}{0}
\setcounter{table}{0}

This section is a brief review of the dimensional splitting SL WENO algorithm for solving the VP system \cite{Qiu_Christlieb, qiu_shu_sl}. The proposed algorithm is up to ninth order accurate in phase space, but subject to dimensional splitting error. The main procedures are outlined below. For simplicity, we consider the problem with only 1-D in space and 1-D in velocity. The domain $(x, v) \in [0, L] \times [-v_{max}, v_{max}]$. The boundary condition is periodic in the $x$-direction, and zero boundary condition is assumed in the $v$-direction.

\bigskip
\noindent
\underline{Dimensional splitting.}  The time splitting form of equation \eqref{eqn:vlasov3d}  is,
\beq
f_t + v\cdot f_x = 0~,~
\label{vlasov_x}
\eeq
\beq
f_t + E(t,{x})\cdot f_v = 0 ~.~
\label{vlasov_v}
\eeq
The splitting form of equation \eqref{eqn:vlasov3d} is first order accurate if one evolves equation \eqref{vlasov_x} for a full time step, then solving the Poisson's equation \eqref{eq: poisson} for electrostatic field $E({ x})$ and finally evolves equation \eqref{vlasov_v} for another full time step in a sequential way. It can be designed to be second order accurate in time by solving equation \eqref{vlasov_x} for a half time step, the solving the Poisson's equation for $E({ x})$, then solving equation \eqref{vlasov_v} for a full time step, followed by  solving equation \eqref{vlasov_x} for a second half time step, i.e., Strang splitting \cite{strang1968construction}.

\bigskip
\noindent
\underline{SL WENO scheme.} The splitting strategy decouples the nonlinearity of the Vlasov equation into two linear hyperbolic equations \eqref{vlasov_x} and \eqref{vlasov_v}. The linearity of 1-D equations allows for a simple implementation of high order SL methods. For example, the conservative SL finite difference WENO schemes proposed for 1-D equations such as equation \eqref{vlasov_x} or \eqref{vlasov_v} have been shown to be successful due to the conservative nature of the scheme formulation, the flexibility of working with point values, the stability and robustness of WENO reconstruction and the large time step efficiency of SL methods \cite{qiu_shu_sl}. We briefly review the methodology below and refer readers to \cite{qiu_shu_sl} for more details.

We adopt the following notations for the numerical discretization. The phase space is discretized by the following uniformly distributed grid points
$$0 = x_0 <  \cdots< x_i \cdots< x_{N_x} = L,$$
$$-v_{max} = v_0 < \cdots< v_j< \cdots < v_{N_v} = v_{max},$$
where $x_i = i \cdot \Delta x$ with mesh size $\Delta x = L/N_x$ and $v_j =-v_{max} + j \cdot \Delta v$ with mesh size $\Delta v = 2v_{max}/N_v$.
Let $x_{i+1/2} = (x_{i}+x_{i+1})/2$, $v_{j+1/2} = (v_{j}+v_{j+1})/2$. Let $f^n_{ij}$ be the point value of numerical solution at $x=x_i$, $v=v_j$ and $t = n \Delta t := t_n$.

Take equation \eqref{vlasov_x} for example: the SL finite difference WENO scheme is based on integrating the conservative form of equation \eqref{vlasov_x} for some fixed $v_j$, over $[t_n, t_{n+1}]$,
\[
f(t_{n+1},x, v_j) = f(t_n, x, v_j ) - \left( \int_{t_n}^{t_{n+1}} v_j f(\tau ,x, v_j ) d\tau \right)_x.
\]
For simplicity of notation, we omit the $v_j$-dependence below under the dimensional splitting framework.
Let $f^n_i \doteq f(t^{n},x_i, v_j)$.
Evaluating the above equation at the grid point $x_i$ gives
\begin{align}
\label{eq: slfd1}
f^{n+1}_i & = f^n_i -  \left.\left( \int_{t_n}^{t_{n+1}} v f(\tau ,x )d\tau \right)_x\right\vert_{x=x_i} = f^n_i -  \mathcal{F}_x\Big|_{x=x_i},
\end{align}
where $\mathcal{F}(x) \doteq \int_{t_n}^{t_{n+1}} vf(\tau ,x )d\tau$.
Let $\mathcal{H}(x)$ be a function whose sliding average is $\mathcal{F}(x)$, i.e.,
\beq
\label{eq: H}
\mathcal{F}(x) = \frac{1}{\Delta x} \int_{x-\frac{\Delta x}{2}}^{x+\frac{\Delta x}{2}} \mathcal{H}(\xi) d\xi.
\eeq
Taking the $x$ derivative of the above equation gives
\[
\mathcal{F}_x = \frac{1}{\Delta x} \left( \mathcal{H}(x+\frac{\Delta x}2) -
\mathcal{H}(x-\frac{\Delta x}2) \right) .
\]
Therefore the equation (\ref{eq: slfd1}) can be written in a conservative form as
\beq
f^{n+1}_i = f^n_i - \frac{1}{\Delta x}( \mathcal{H}(x_{i+\frac12}) -  \mathcal{H}(x_{i-\frac12})),
\eeq
where $ \mathcal{H}(x_{i+\frac12})$ is called the flux function. By following information along characteristics,
$
\mathcal{F}(x_i) = \int_{t_n}^{t_{n+1}} vf(\tau ,x_i )d\tau = \int_{x_i^\star}^{x_i} f d\xi
$
can be evaluated by reconstructing function $f$ at $t_n$ from neighboring point values. We denote this reconstruction procedure as $\mathcal{R}_1$.
Similar to the idea in the finite difference WENO scheme, $\mathcal{H}(x_{i+\frac12})$ can be reconstructed in high order from several of its neighboring cell averages
\[
\bar{\mathcal{H}}_k = \frac{1}{\Delta x} \int_{x_{k-\frac12}}^{x_{k+\frac12}} \mathcal{H}(\xi) d\xi \stackrel{\eqref{eq: H}}{=} \mathcal{F}(x_k), \quad
k=i-p, \cdots, i+q.
\]
We denote this reconstruction procedure as $\mathcal{R}_2$.
In summary, a SL finite difference scheme in evolving equation \eqref{vlasov_x} from $t_n$ to $t_{n+1}$ can be designed as follows:
\begin{enumerate}
\item \vspace{-.1in}
At each of the grid points at time level $t_{n+1}$, say $(x_{i}, t_{n+1})$,  trace the characteristic back to time level $t_n$ at $x^\star_{i} = x_i - v \Delta t$.
\item \vspace{-.1in}
Reconstruct $\mathcal{F}(x_i) = \int_{t_n}^{t_{n+1}} vf(\tau ,x_i )d\tau = \int_{x_i^\star}^{x_i} f d\xi$ from $\{f^n_i\}_{i=1}^{N_x}$. We use $\mathcal{R}_1$
to denote this reconstruction procedure
$\mathcal{R}_1[x_i^\star, x_i] (f^n_{i-p_1}, \cdots, f^n_{i+q_1})
$
in approximating $\mathcal{F}(x_i)$, where $(i-p_1, \cdots, i+q_1)$ indicates the stencil used in the reconstruction. $\mathcal{R}_1[a, b]$ indicates the reconstruction of $\int_{a}^b f(t, \xi)d\xi$.
\item \vspace{-.1in}
Reconstruct $\{\mathcal{H}(x_{i+\frac12})\}_{i=0}^{N_x}$ from $\{\bar{\mathcal{H}}_i\}_{i=1}^{N_x}$. We use $\mathcal{R}_2$
to denote this reconstruction procedure
$
\hat{F}_{i+\frac12}\doteq \mathcal{R}_2 (\bar{\mathcal{H}}_{i-p_2}, \cdots, \bar{\mathcal{H}}_{i+q_2})
$
in approximating $\mathcal{H}(x_{i+\frac12})$, where $(i-p_2, \cdots, i+q_2)$ indicates the stencil used in the reconstruction.
\item \vspace{-.1in}
Update the solution $\{f^{n+1}_i\}_{i=1}^{N_x}$ by
\beq
\label{eq: 1d_cons_2}
f^{n+1}_i = f^n_i - \frac{1}{\Delta x}( \hat{F}_{i+\frac12} - \hat{F}_{i-\frac12}),
\eeq
with numerical fluxes $\hat{F}_{i\pm\frac12}$ computed in the previous step.
\end{enumerate}
When the reconstruction stencils in $\mathcal{R}_1$ and $\mathcal{R}_2$ above only involve one neighboring point value of the solution, then the scheme is first order
accurate in space. In fact, the scheme reduces to a first order upwind scheme when the time step is within CFL restriction. The proposed SL finite difference
scheme can be designed to be of high order accuracy by including more points in the  stencil for $\mathcal{R}_2 \circ \mathcal{R}_1$, the composition of $\mathcal{R}_1$ and $\mathcal{R}_2$, to reconstruct the numerical flux
\begin{equation}
\label{weno_flux}
\hat{F}_{i+\frac12} = \mathcal{R}_2 \circ \mathcal{R}_1 (f^n_{i-p}, \cdots, f^n_{i+q}),
\end{equation}
where $(i-p, \cdots, i+q)$ indicates the stencil used in the reconstruction process.
The WENO mechanism can be introduced
in reconstruction procedures as a stable and non-oscillatory method to capture fine scale structures. The details of the design of high order SL WENO algorithm are described in \cite{qiu_shu_sl}. It is numerically demonstrated in \cite{Qiu_Christlieb, qiu_shu_sl, Qiu_Shu2} that the proposed high order SL WENO method works very well in Vlasov simulations with extra large time step evolution.

The computational procedure for the first order dimensional splitting SL finite difference WENO scheme for the VP system is summarized by the flow chart below.

\bigskip
\fbox{
\begin{minipage}[htb]{0.9\linewidth}
\textbf{Algorithm 1:}
A first order dimensional splitting SL finite difference WENO scheme for the VP system:
\begin{enumerate}
\item Evolve the solution $f^n$ by solving equation \eqref{vlasov_x} to advect in $x$-direction for a time step $\Delta t$ with a SL finite difference WENO scheme \eqref{eq: 1d_cons_2}, where
the numerical flux $\hat{F}_{i+\frac12}$ is reconstructed by \eqref{weno_flux} in a WENO fashion. Denote the solution after the evolution by $f^{n,*}$.
\item Solve electrostatic field $E^{n,*}$ induced by $f^{n,*}$ from the Poisson's equation using a fast Fourier transform (FFT).
\item Evolve the solution $f^{n,*}$ by solving equation \eqref{vlasov_v} to advect in $v$-direction for a time step $\Delta t$ to get $f^{n+1}$ as step {1}.
\end{enumerate}
\end{minipage}
}

\bigskip
\begin{rem}
\textbf{Algorithm 1} can be improved into second order in time by adopting Strang splitting as we mentioned at the beginning of Section 2.
\end{rem}

\begin{rem}
There are different conservative SL finite difference procedures proposed in  \cite{Qiu_Christlieb, qiu_shu_sl, Qiu_Shu2}. These approaches have similar performance for linear advection equations with constant coefficients. However, we choose to review the SL procedure from \cite{qiu_shu_sl},
as it is more general in the sense that it can be applied to linear advection equations with variable coefficients. For solving the VP system, the methodologies proposed in \cite{Qiu_Christlieb, qiu_shu_sl, Qiu_Shu2} can be applied. However, for the guiding center model in Section~\ref{sec: 3.3}, the dimensional splitting equations have variable coefficients. Hence, the above reviewed procedure from \cite{qiu_shu_sl} is needed.
\end{rem}

\section{IDC methods for improving the temporal order of accuracy}
\label{sec3}
\setcounter{equation}{0}
\setcounter{figure}{0}
\setcounter{table}{0}

In this section, we first briefly review the IDC procedure for initial value problems \cite{DGR, highRKsdc}. Then we introduce the proposed IDC scheme in reducing the dimensional splitting error for the VP system in Section~\ref{sec: 3.2} and for the guiding center model in Section ~\ref{sec: 3.3}.

\subsection{Overview of IDC methods}
\label{sec3.1}

We provide a brief review of IDC methods \cite{highRKsdc}, designed for scalar/system of initial value problems in the following form,
\begin{equation}
\label{ode}
\left\{\begin{array} {l}
\displaystyle
y'(t)  = g(t, y), \quad t \in[0, T], \\
\displaystyle
y(0) =y_0.
         \end{array}
   \right.
\end{equation}
The time domain, $[0, T]$, is discretized into intervals,
\begin{equation*}
\label{mesh1}
0 = t_1 < t_2 < \cdots < t_n < \cdots <t_N = T,
\end{equation*}
and each interval, $I_n = [t_n, t_{n+1}]$, is further discretized into
sub-intervals,
\begin{equation}
\label{eq: mesh2}
t_n = t_{n,0} < t_{n, 1} < \cdots < t_{n,m} < \cdots < t_{n, M} = t_{n+1}.
\end{equation}
The IDC method on each time interval $[t_n, t_{n+1}]$ is described
below. We drop the subscript $n$, e.g., $\tau_{0} : = t_{n, 0}$ in
\eqref{eq: mesh2}, with the understanding that the IDC method is described
for one time interval. We also refer to $\tau_m := t_{n,m}$ as grid
points or quadrature nodes, whose index $m$ runs from $0$ to $M$.
In the IDC method, the size of sub-intervals are uniform.
Let $\Delta \tau \doteq \frac{t_{n+1}-t_n}{M}$, then $\tau_m = t_n + m \Delta \tau$, $m=0, \ldots, M$.
The procedure of an IDC method with $M+1$ uniformly distributed quadrature nodes as in equation \eqref{eq: mesh2}
and with $K$ correction loops denoted as IDC$(M+1)$J$(K)$ is the following:
\begin{enumerate}
\item (prediction step) Use a low order numerical method to
  obtain a numerical solution, $\vec{\eta}^{[0]} = (\eta^{[0]}_{0},
   \ldots,\eta^{[0]}_{m},\ldots, \eta^{[0]}_{M}), $
  which is a low order approximation to the exact solution at quadrature points.
   For example, applying a first order forward
  Euler method to \eqref{ode} gives
  \[
  \eta^{[0]}_{m+1} =
  \eta^{[0]}_{m} + \Delta \tau g(t, \eta^{[0]}_m),\, m=0,\ldots,M-1.
  \]
\item (correction loops) Use the error function to improve the accuracy
  of the scheme at each iteration.\\
  For $k = 1,\ldots, K$ ($K$ is the number of correction steps)
\begin{enumerate}
\item Denote {\em the error function} from the previous step as
  \begin{align*}
    e\kfm(t)  = y(t) - \eta\kfm(t),
    \label{err_func}
  \end{align*}
  where $y(t)$ is the exact solution and $\eta\kfm(t)$ is an $M^{th}$
  degree polynomial interpolating $\vec{\eta}\kvm$.
  \item Denote {\em the residual function} as
\begin{equation*}\epsilon\kfm(t) =
  (\eta\kfm)'(t) - g(t, \eta\kfm(t)).
\end{equation*}
\item Compute {\em the numerical error vector}, $\vec{\delta}\kv =
  (\delta\kv_0, \ldots, \delta\kv_m, \ldots, \delta\kv_M)$, using a
  low order numerical method to discretize the following {\em error equation} with a zero initial condition,
  \begin{align}
  \label{ode_err}
  \left(e^{(k-1)} + \int_0^t\!\!\epsilon^{(k-1)}(\tau)\,d\tau \right)'\!\!(t)
    &=  g(t, \eta\kfm(t) + e\kfm(t)) -
      g(t, \eta\kfm(t)).
  \end{align}
For example, applying a first order forward Euler scheme to equation \eqref{ode_err} gives
\begin{align}
  \label{eqn:error_sdcfe}
  \delta\kv_{m+1}
  =\delta\kv_m + \Delta \tau (g(\tau_m, \eta\kvm_m+\delta\kv_m) -
  g(\tau_m, \eta_m\kvm)) + \sum_{\ell=0}^M{\alpha_{m,\ell}~ g(\tau_\ell,\eta_\ell^{[k-1]})}  \nonumber\\
  +\eta_m^{[k-1]} - \eta_{m+1}^{[k-1]}, \quad
   m = 0,\ldots,M-1.
\end{align}
To get equation \eqref{eqn:error_sdcfe} from the discretization of \eqref{ode_err}, $\int_{\tau_m}^{\tau_{m+1}} \epsilon\kfm(t)\,dt$ is approximated by
\[
\eta_{m+1}^{[k-1]}-\eta_m^{[k-1]}-  \sum_{\ell=0}^M{\alpha_{m,\ell} ~g(\tau_\ell,\eta_\ell^{[k-1]})},
\]
where $\sum_{j=0}^M{\alpha_{m,\ell}~ g(\tau_\ell,\eta_\ell^{[k-1]})}$ approximates $\int_{\tau_m}^{\tau_{m+1}} g(t, \eta\kfm(t)) dt$ by quadrature formulas. We note that such a way of evolving the error function in the IDC algorithm is advantageous compared with that in the traditional DC algorithm. It introduces more stability by using the integral form, rather than the differentiation form, of the residual.
\item Update the numerical solution $\vec{\eta}\kv = \vec{\eta}\kvm +  \vec{\delta}\kv$.
\end{enumerate}
\end{enumerate}
\begin{rem}  (About notations.) In our description of IDC, we let English letters $y_m$, $e^{(k)}_m$ denote the exact solutions and exact error functions, and Greek letters ${\eta}^{[k]}_m$, ${\delta}^{[k]}_m$ denote the numerical approximations to the exact solutions and error functions. We use subscript $m$ to denote the location $t = \tau_m$. The superscripts $k$ with round brackets ($(k)$) and square brackets ($[k]$) are for functions and vectors (or their components) respectively at the prediction ($k=0$) and correction loops ($k=1, \ldots, K$). We let $\vec{\cdot}$ denote the vector on IDC quadrature nodes. For example, $\vec{y} = (y_0, \ldots, y_M)$.
\end{rem}

\begin{rem} (About the distribution of the quadrature nodes.) The IDC methods reviewed above adopt the uniformly distributed quadrature nodes to compute the residual. In \cite{highRKsdc}, it is proved under some mild assumption, that the order accuracy of a IDC method can be increased by $r$ order when an $r^{th}$ order Runge-Kutta integrator is used to solve the error equation in each correction loop. The numerical results reported in \cite{christlieb2009comments} show that such high order accuracy lifting property does not always hold for an SDC method, which is constructed with Gaussian quadrature nodes.
\end{rem}

\begin{rem}(About computational cost and storage requirement.)
In terms of the order of accuracy, per $\Delta \tau$ (a subinterval size in the IDC method),
numerical solutions $\eta^{[K]}_m$ are $(K+1)^{th}$ order approximations to the exact solution at quadrature nodes $y_m$ for $m=0, \ldots, M$; in terms of computational cost, per $\Delta \tau$, there are $(K+1)$ function evaluations for a $(K+1)^{th}$ order IDC$(M+1)$J$(K)$ method. In this sense, the IDC method
is considered to be efficient with relatively low computational cost among Runge-Kutta methods with the same order of accuracy. We remark that an IDC method can be considered as a one-step Runge-Kutta method with its Butcher table constructed in \cite{christlieb2009comments}. At the same time, an IDC$(M+1)$J$(K)$ method requires storage space for numerical solutions at $(M+1)$ quadrature nodes.
\label{rem: cost}
\end{rem}

\subsection{IDC methods for the VP system}
\label{sec: 3.2}

The dimensional split SL method described in Section~\ref{sec2} is very high order accurate in phase spaces, but only low order in time. In this subsection, we use the IDC framework to increase the temporal order of the dimensional splitting method. We will only consider using first order splitting for brevity, but also comment on the use of Strang splitting. Consider the VP system with only 1-D in space and 1-D in velocity for simplicity of notation.
\begin{enumerate}
\item (prediction step) Use the proposed dimensional splitting SL WENO method ({\bf Algorithm 1}) described in Section~\ref{sec2} in the prediction step of the IDC framework. More specifically, predict solution
$
\vec{\eta}^{[0]} = (\eta_{0}^{[0]}, \ldots, \eta_{M}^{[0]})
$
at time step subintervals \eqref{eq: mesh2} for each spatial and velocity grid point, say
$
(x_i, v_j), \quad \forall i=1, \ldots, N_x, \quad j=-N_v/2, \ldots, N_v/2.
$
\item (correction loops) Use the error function to reduce the dimensional splitting error at each iteration. Our correction procedure is based on a fixed location, say $(x_i, v_j)$.\\
For $k=1,\ldots, K$ ($K$ is the number of correction steps)
\begin{enumerate}
\item {\em The error function} is defined as $e^{(k-1)}(t, x_i, v_j) = f(t, x_i, v_j) - \eta^{(k-1)}(t, x_i, v_j)$,
where $f(t, x_i, v_j)$ is the exact solution and $\eta^{(k-1)}(t, x_i, v_j)$ is the polynomial interpolating $\vec{\eta}^{[k-1]} = (\eta_{0}^{[k-1]}(x_i, v_j), \ldots, \eta_{M}^{[k-1]}(x_i, v_j))$ at quadrature points \eqref{eq: mesh2} over a time interval $[t_n, t_{n+1}]$.
\item {\em The residual function} is defined as
\beq
\label{eq: residual}
\epsilon(t, x, v) = - (\eta_t + v \cdot \eta_x + {E}^{\eta}({x}) \cdot \eta_v),
\eeq
where ${E}^{\eta}$ is the electrostatic field induced by the numerical distribution function $\eta(t, x, v)$.
\item {\em The error equation} about the error function is obtained by adding the residual equation \eqref{eq: residual} to the Vlasov equation \eqref{eqn:vlasov3d},
\beq
\label{eq: err}
e_t + v \cdot e_x + (E^{\eta}+E^{e}) \cdot e_v+E^{e} \cdot \eta_v = \epsilon,
\eeq
where $E^{e}$ is the electrostatic field induced by the error function $e(t, x, v)$.
\item Evolve the error equation \eqref{eq: err} with zero initial condition by the same dimensional splitting SL WENO method for spatial advection and velocity acceleration/deceleration as that for the Vlasov equation. Specifically, we split the error equation \eqref{eq: err} into three parts,
\beqa
\label{eq: spatial_corr}
e_t + v \cdot e_x = 0 && \mbox{(spatial advection)} \\
\label{eq: veloc_corr}
e_t  + E^{\eta+e} \cdot e_v = 0 && \mbox{(velocity ac/deceleration)} \\
\label{eq: source}
e_t+E^{e} \cdot \eta_v = \epsilon && \mbox{(source term)}.
\eeqa
To evolve the error equation \eqref{eq: err} from $\tau_m$ to $\tau_{m+1}$ in a splitting fashion, we first evolve the solution $\delta^{[k]}_{m}$ by approximating equation \eqref{eq: spatial_corr} with the SL WENO scheme as in {\bf Algorithm 1}. The solution after the evolution is denoted by $\delta^{[k]}_{m, *}$.
After that, we solve the electrostatic field $E^{\eta^{[k-1]}_m+\delta^{[k]}_{m,*}}$ induced by $\eta^{[k-1]}_m+\delta^{[k]}_{m,*}$ from the Poisson's equation, then get $E^{\delta^{[k]}_{m,*}}=E^{\eta^{[k-1]}_m+\delta^{[k]}_{m,*}}-E^{\eta^{[k-1]}_m}$. We then evolve $\delta^{[k]}_{m, *}$ by approximating equation \eqref{eq: veloc_corr} again with the SL WENO scheme as in {\bf Algorithm 1}. The numerical solution is denoted by $\delta^{[k]}_{m,**}$. Finally we solve equation \eqref{eq: source}, but using the integral form of the residual in an SDC/IDC fashion,
  \begin{align}
  \label{eq: idc_cor_sour_vp}
    \left(e(t, x, v) - \int_{t_n}^t\!\!\epsilon(\tau, x, v)\,d\tau \right)'\!\!(t)
    &= - E^{e} \cdot \eta_v.
  \end{align}
Note that $E^e \cdot \eta_v$, as well as the terms $v \cdot \eta_x$ and ${E}^\eta \cdot \eta_v$ in equation \eqref{eq: residual} for the residual $\epsilon$, are approximated in a flux difference form to ensure mass conservation. In the simulation, we use a fifth order WENO procedure to reconstruct all the numerical fluxes.
Similar to \eqref{eqn:error_sdcfe}, we approximate \eqref{eq: idc_cor_sour_vp} at $(x_i, v_j)$ by
\begin{align}
  \label{eqn:error_sdcfe_vp}
  \delta\kv_{m+1}
  =\delta^{[k]}_{m, **} - \frac{\Delta\tau}{\Delta v}E^{\delta^{[k]}_{m,*}}(\widehat{{\eta}_m^{[k-1]}}_{i,j+\frac12} - \widehat{{\eta}_m^{[k-1]}}_{i,j-\frac12})
 - \sum_{\ell=0}^M{\alpha_{m,\ell}~ g(\eta_\ell^{[k-1]})} \notag\\
  +\eta_{m+1}^{[k-1]}-\eta_m^{[k-1]} , \quad
   m = 0,\ldots,M-1,
\end{align}
with
\begin{equation}
\label{eq: g}
g(\eta_\ell^{[k-1]})=v \frac{\widehat{\eta_\ell^{[k-1]}}_{i+\frac12,j}-\widehat{\eta_\ell^{[k-1]}}_{i-\frac12,j}}{\Delta x}+{E}^{\eta^{[k-1]}_\ell} \frac{\widehat{\eta_\ell^{[k-1]}}_{i,j+\frac12}- \widehat{\eta_\ell^{[k-1]}}_{i,j-\frac12}}{\Delta v},
\end{equation}
where we omit the $(i, j)$ dependence when there is no confusion.
\item Update the solution by the approximate error function computed from the correction step,
\begin{equation}
\label{eq: update}
\vec{\eta}\kv(x_i, v_j) = \vec{\eta}\kvm(x_i, v_j) +  \vec{\delta}\kv(x_i, v_j).
\end{equation}
\end{enumerate}
\end{enumerate}

The flow chart of the SL WENO algorithm combined with IDC$(M+1)$J$(K)$ for the VP system is outlined below.

\bigskip
\fbox{
\begin{minipage}[htb]{0.9\linewidth}
\textbf{Algorithm 2:} A SL finite difference WENO scheme coupled in the IDC framework:
\begin{enumerate}
\item Find the prediction solution $\vec{\eta}^{[0]}=(\eta^{[0]}_0,\ldots,\eta^{[0]}_M)$ at time step subintervals for each spatial and velocity location $(x_i, v_j),~\forall i=1, \ldots, N_x,~ j=-N_v/2, \ldots, N_v/2$, by using the dimensional splitting SL WENO scheme as in {\bf Algorithm 1}.
\item {\bf For} $k=1,\ldots,K$

    {\em Perform the correction loop to update the solution $\vec{\eta}\kv$.}
    \begin{itemize}
    \item Solve the numerical error vector $\vec{\delta}\kv=(\delta^{[k]}_0,\ldots,\delta^{[k]}_M)$ by evolving the split error equations \eqref{eq: spatial_corr}-\eqref{eq: source} with zero initial condition at time step subintervals for each spatial and velocity location. Specifically,

{\bf For} $m=0,\ldots,M-1$,
\begin{enumerate}
\item Evolve $\delta\kv_{m}$ by solving equation \eqref{eq: spatial_corr} to get $\delta\kv_{m,*}$ by the SL WENO as in {\bf Algorithm 1}.
\item Evolve $\delta\kv_{m, *}$ by solving equation \eqref{eq: veloc_corr} to get $\delta\kv_{m,**}$ by the SL WENO as in {\bf Algorithm 1}.
\item Evolve $\delta\kv_{m,**}$ by solving equation \eqref{eq: source} in an IDC fashion using equation \eqref{eqn:error_sdcfe_vp} to get $\delta\kv_{m+1}$.
\end{enumerate}
{\bf End For}
\item update the solution by $\vec{\eta}\kv = \vec{\eta}\kvm +  \vec{\delta}\kv$.
\end{itemize}
{\bf End For}
\end{enumerate}
\end{minipage}
}

\bigskip

The proposed IDC - dimensional splitting SL WENO scheme enjoys the mass conservation property; see the following proposition.
\begin{prop}
\label{prop: mass}
{\bf Algorithm 2} conserves the total mass for solving the VP system, if the boundary conditions are periodic.
\end{prop}
\noindent Proof:
First note that {\bf Algorithm 1} conserves the total mass when periodic boundary conditions are imposed.
Thus, we have
\begin{align}
\label{eq: proof_prediction}
\sum_i\sum_j\eta^{[0]}_0(x_i,v_j)=\sum_i\sum_j\eta^{[0]}_1(x_i,v_j)=\ldots = \sum_i\sum_j\eta^{[0]}_M(x_i,v_j),
\end{align}
for the prediction. By \eqref{eq: update}, to prove the mass conservation, it is sufficient to prove
\begin{align}
\label{eq: proof_cor}
\sum_i\sum_j\delta^{[k]}_m(x_i,v_j) = 0, \quad k=1,\ldots,K, \quad m = 0,\ldots,M.
\end{align}

Let's prove \eqref{eq: proof_cor} for $k=1$. A similar argument carries over for general $k$.
The split error equations \eqref{eq: spatial_corr}-\eqref{eq: source} are solved with zero initial condition, hence
$
\sum_i\sum_j\delta^{[1]}_0(x_i,v_j) = 0.
$
Now assume
$
\sum_i\sum_j\delta^{[1]}_m(x_i,v_j) = 0,
$
we will prove
$
\sum_i\sum_j\delta^{[1]}_{m+1}(x_i,v_j) = 0.
$
Due to total mass conservative property of {\bf Algorithm 1}, we get
\begin{align}
\label{eq: proof_med}
\sum_i\sum_j\delta^{[1]}_{m,*}(x_i,v_j) = 0, \quad \sum_i\sum_j\delta^{[1]}_{m,**}(x_i,v_j) = 0.
\end{align}
Since the derivatives in equation \eqref{eqn:error_sdcfe_vp} are written into a flux difference form,
it follows that
\begin{align*}
  \sum_i\sum_j\delta^{[1]}_{m+1}(x_i,v_j)
  =\sum_i\sum_j\delta^{[1]}_{m,**} (x_i,v_j)-\sum_i\sum_j
   \frac{\Delta\tau}{\Delta v}E^{\delta^{[k]}_{m,*}}(\widehat{{\eta}_m^{[k-1]}}_{i,j+\frac12} - \widehat{{\eta}_m^{[k-1]}}_{i,j-\frac12})
  \nonumber\\
  - \sum_i\sum_j\sum_{\ell=0}^M{\alpha_{m,\ell}~ g(\eta_\ell^{[0]}(x_i,v_j))}+\sum_i\sum_j\eta_{m+1}^{[0]}(x_i,v_j)-\sum_i\sum_j\eta_m^{[0]}(x_i,v_j).
\end{align*}
Thanks to \eqref{eq: proof_prediction}, \eqref{eq: proof_med} and the cancellation of the flux difference form, we get
\begin{align*}
\sum_i\sum_j\delta^{[1]}_{m+1}(x_i,v_j) &= - \sum_i\sum_j\sum_{\ell=0}^M{\alpha_{m,\ell}~ g(\eta_\ell^{[0]}(x_i,v_j))}\\
& = -\sum_{\ell=0}^M{\alpha_{m,\ell}~\sum_i\sum_j g(\eta_\ell^{[0]}(x_i,v_j))}\\
& = 0,
\end{align*}
where the last equality holds due to the flux difference form of $g$ in \eqref{eq: g}. By induction, we complete the proof. $\Box$

\begin{rem}
\textbf{Algorithm 2} can be extended to IDC methods coupled with second order Strang splitting without additional complication. The only modification is to
employ the Strang splitting SL finite difference WENO scheme to get a prediction, and again employ the Strang splitting to solve the split error equation \eqref{eq: spatial_corr}-\eqref{eq: source}.
The procedure of an IDC method coupled with the second order Strang splitting is denoted as IDC-Strang$(M+1)$J$(K)$, when $M+1$ uniformly distributed quadrature nodes
and $K$ correction loops are used. The numerical results reported in Section 5 indicate that the temporal order accuracy of the IDC-Strang method for the VP system can be increased with second order per correction. The modified IDC-Strang splitting SL WENO scheme also enjoys the mass conservation property. The proof is quite similar, therefore we omit it.
\end{rem}

\subsection{IDC methods for the guiding center model}
\label{sec: 3.3}

We consider the guiding center model, which describes a highly magnetized plasma in the transverse plane of a tokamak \cite{sonnendruecker,crouseilles2010conservative}. We consider equation
\begin{equation}
\label{guiding}
\rho_t + \mathbf{E}^{\perp}\cdot\nabla\rho = 0.
\end{equation}
where $\rho$ is the particle density function, $\mathbf{E}^{\perp}=(-E_2,E_1)$ with the electrostatic field $\mathbf{E}=(E_1,E_2)$ satisfying a Poisson's equation
\begin{equation}
\Delta \Phi =\rho,\quad \mathbf{E}=-\nabla\Phi.
\end{equation}
Compared to the VP system, the 1-D equations obtained from dimensional splitting of equation \eqref{guiding} are variable coefficient equations. We apply the SL finite difference WENO algorithm proposed in \cite{qiu_shu_sl} as described in the previous section as 1-D solvers.
The 2-D Poisson's equation is again solved by a 2-D FFT. The computational procedure of applying the IDC method to reduce the splitting error is similar to that for the VP system, except that we need to formulate a new residual function and an error equation for the guiding center model \eqref{guiding}. Below we
provide the residual function and the error equation for the guiding center Vlasov equation with notations that are consistent with the previous subsection.
\bit
\item
{\em The residual function} is defined as
\beq
\label{residual}
\epsilon(t, x, y) = - \left(\eta_t - (E_2^\eta \eta)_x + (E_1^\eta \eta)_y\right)
\eeq
where $\mathbf{E}^\eta = -\nabla \Phi^\eta$ with
$
\Delta\Phi^\eta = \eta.
$

\item {\em The error equation} about {\em the error function} $e(t, x, y) = \rho(t, x, y) - \eta(t, x, y)$ is obtained by adding \eqref{residual} to \eqref{guiding},
\beq
e_t - (E^{\rho}_2  \rho - E^\eta_2  \eta)_x + (E^{\rho}_1  \rho - E^\eta_1  \eta)_y  =  \epsilon,
\eeq
where $\mathbf{E}^{\rho}$ and $\mathbf{E}^{\eta}$ are the electrostatic field induced by the exact solution $\rho(t, x, y)$ and the numerical solution $\eta(t, x, y)$ respectively.
From $\rho = \eta + e$ and $\mathbf{E}^{\rho} = \mathbf{E}^{\eta} + \mathbf{E}^e$,  we have
\beq
e_t - \left((E^e_2+E^\eta_2)  e\right)_x + \left((E^e_1+E^\eta_1)  e\right)_y - (E^e_2  \eta)_x  + (E^e_1 \eta)_y =  \epsilon,
\eeq
where $\mathbf{E}^e$ is the electrostatic field induced by the error function $e(t, x, y)$. Similar to the proposed scheme for the VP system, the error equation is evolved with zero initial conditions by dimensional splitting,
\beqa
\label{eq: spatial_corr_guiding}
e_t - \left((E^e_2+E^\eta_2)  e\right)_x  &=& 0, \\
\label{eq: veloc_corr_guiding}
e_t + \left((E^e_1+E^\eta_1)  e\right)_y &=& 0, \\
\label{eq: source_guiding}
e_t - (E^e_2  \eta)_x  + (E^e_1  \eta)_y &=&  \epsilon.
\eeqa
Again, the 1-D equations \eqref{eq: spatial_corr_guiding}-\eqref{eq: veloc_corr_guiding} are solved by the SL WENO method as in {\bf Algorithm 1}, and the last equation above is solved in a similar manner as equation \eqref{eq: idc_cor_sour_vp} for the VP system. Similar to Proposition~\ref{prop: mass}, the algorithm enjoys the mass conservation property, since the splitting is performed in a conservative way.
\eit

\section{Numerical stability}
\label{sec4}
\setcounter{equation}{0}
\setcounter{figure}{0}
\setcounter{table}{0}

IDC is a numerical approach in generating time stepping algorithms with high order accuracy, yet numerical stability of the IDC method using the SL WENO as the base scheme remains to be investigated. In the following, we investigate stability properties of the proposed IDC SL WENO method via classical Fourier analysis. We provide the CFL restriction for stability when the method is applied to a linear problem as guidance for choosing numerical time steps for general nonlinear problems in Section~\ref{sec5}.

We consider a linear model problem:
\begin{equation}
\label{eq: model_cfl}f_t+ f_x=0.
\end{equation}
Assume the mesh is uniform and boundary condition is periodic. We consider a subinterval in IDC with time step size $\Delta \tau$.
An explicit linear scheme for equation \eqref{eq: model_cfl} can be written in the following form:
\begin{equation}
\label{eq: scheme_cfl}
f_j^{n+1} = \sum_{k=-r}^l C_k f_{j+k}^n,
\end{equation}
where $C_k,\,k=-r,\ldots,l$ are constants, that depend on the CFL number $\lambda \doteq \Delta \tau/\Delta x$ but are independent of the solution. For example, the third order linear SL scheme combined with IDC2J0 reads:
\begin{align*}
\displaystyle f_j^{n+1} = & \frac{1}{6} \lambda^\star \left((\lambda^\star)^2-1\right) f^n_{j^\star-2}+\frac{1}{2} \lambda^\star \left(-(\lambda^\star)^2+\lambda^\star+2\right) f^n_{j^\star-1} \\[2mm]
\displaystyle& +\frac{1}{2} \left((\lambda^\star)^3-2 (\lambda^\star)^2-\lambda^\star+2\right) f^n_{j^\star} -\frac{1}{6} \lambda^\star \left((\lambda^\star)^2-3 \lambda^\star+2\right) f^n_{j^\star+1}
\end{align*}
where $j^\star = j-\lfloor\lambda\rfloor$ and $\lambda^\star=\lambda-\lfloor\lambda\rfloor$. The third order linear SL scheme combined with IDC2J1 reads:
\begin{align*}
f_j^{n+1} = & -\frac{1}{72} \lambda^2 \left(\lambda^2-1\right) f^n_{j-4} +  \frac{1}{24} \lambda^2 \left(3 \lambda^2-\lambda-4\right) f^n_{j-3}\\[2mm]
& +\frac{1}{12} \lambda \left(-4 \lambda^3+4 \lambda^2+7 \lambda-2\right) f^n_{j-2} + \frac{1}{36} \lambda \left(13 \lambda^3-24 \lambda^2-16 \lambda+36\right) f^n_{j-1} \\[2mm]
& -\frac{1}{24} \left(3 \lambda^4-10 \lambda^3+5 \lambda^2+12 \lambda-24\right) f^n_j + \frac{1}{24} \lambda \left(-\lambda^3+\lambda^2+4 \lambda-8\right) f^n_{j+1} \\[2mm]
& + \frac{1}{36} \lambda^2 \left(\lambda^2-3 \lambda+2\right) f^n_{j+2}
\end{align*}
when $0<\lambda<1$; and
\begin{align}
\notag
\displaystyle f_j^{n+1} =& -\frac{1}{72} \lambda^2 \left(\lambda^2-3 \lambda+2\right) f^n_{j-5} + \frac{1}{24} \lambda^2 \left(3 \lambda^2-10 \lambda+7\right) f^n_{j-4}\\[2mm]
\displaystyle\notag & -\frac{1}{12} \lambda^2 \left(4 \lambda^2-16 \lambda+13\right) f^n_{j-3} + \frac{1}{36} \lambda \left(13 \lambda^3-63 \lambda^2+71 \lambda-6\right) f^n_{j-2}\\[2mm]
\displaystyle\notag & +\frac{1}{24} \lambda \left(-3 \lambda^3+19 \lambda^2-34 \lambda+24\right) f^n_{j-1} -\frac{1}{24} \left(\lambda^4-4 \lambda^3+\lambda^2+12 \lambda-24\right) f^n_{j} \\[2mm]
\displaystyle\notag & +\frac{1}{36} \lambda \left(\lambda^3-6 \lambda^2+11 \lambda-12\right) f^n_{j+1}
\end{align}
when $1\leq\lambda<2$.

Classical Fourier analysis is performed to linear schemes in the form of \eqref{eq: scheme_cfl}. We substitute $f_j^n=a^n e^{Ij\xi}$ (where $I=\sqrt{-1}$) into linear schemes in the form of \eqref{eq: scheme_cfl} and compute the corresponding amplification factors $a$. We use Mathematica to derive the explicit form of the schemes and the corresponding amplification factors $a$. The time step restriction from linear stability can be obtained by maximizing the CFL number $\lambda$ with the constraint that $a\le 1$ for any $\xi\in[0,2\pi]$.
We remark that it is exceedingly tedious and sometimes very difficult to derive the upper bounds of the CFL number analytically, especially for high order schemes. Hence we rely on numerical approaches to obtain such upper bounds: we find the maximum $\lambda$ such that $a\leq 1$ for $2000$ evenly distributed points $\{\xi_n\}$ over $[0,2\pi]$. In Table \ref{tab_cfl}, we list the upper bounds of the CFL numbers for several numerical schemes which use linear SL schemes as base solvers in the IDC framework, of different orders. It is observed that the CFL upper bounds for an IDC subinterval are comparable to those of the Eulerian Runge-Kutta WENO scheme with the same orders of accuracy. As commented in Remark~\ref{rem: cost}, there are $(K+1)$ function evaluations for the IDC$(M+1)$J$(K)$ method per IDC subinterval $\Delta \tau$, where $(K+1)^{th}$ order accuracy is achieved at all quadrature nodes $\tau_m$, $m=0, \cdots, M$. The IDC method
is considered to be efficient among Eulerian Runge-Kutta methods with the same order of accuracy. We remark that the Fourier analysis can be extended to arbitrary order cases, however the algebraic manipulations may become prohibitively complicated.

\begin{table} [htp]
\begin{center}
\caption{The upper bounds of the CFL numbers. SL3 and SL5 are SL schemes with the third and fifth order linear reconstructions. IDC$(M+1)$J$(K)$ denotes an IDC procedure with $M+1$ uniformly distributed quadrature nodes and $K$ correction loops.}
\bigskip
\begin{tabular}{|c|c |c|c|c|c |
}\hline
Scheme& IDC2J0 &IDC2J1 &IDC3J0&IDC3J1&IDC3J2\\
\hline
 SL3 & No restriction& 1.50& No restriction& 0.73 &0.67 \\ \hline
 SL5 & No restriction &1.55& No restriction&0.71& 0.66\\
\hline
\end{tabular}
\label{tab_cfl}
\end{center}
\end{table}

\section{
Numerical examples
}
\label{sec5}
\setcounter{equation}{0}
\setcounter{figure}{0}

In this section, we present some simulation results for solving the VP system, the guiding center Vlasov model and the 2-D incompressible Euler equations. Through these examples, we numerically demonstrate the low order dimensional splitting error, and the IDC's ability to correct these errors.
In the simulations,
\[
\text{CFL}  = \Delta \tau \left(\frac{|c_x|}{\Delta x} + \frac{|c_v| }{\Delta v}\right),
\]
where $|c_x|$ and $|c_v|$ are maximum wave propagation speeds in the $x-$ and the $v-$ directions respectively, and $\Delta \tau$ is the size of a sub-interval in the IDC method. Again, IDC$(M+1)$J$(K)$ and IDC-Strang$(M+1)$J$(K)$  denote  IDC procedures with $M+1$ uniformly distributed quadrature nodes and $K$ correction loops coupling the first order splitting and the Strang splitting, respectively.

\subsection{The VP system}
We consider solving the VP system by the proposed scheme
with the following three initial conditions.
\begin{itemize}
\item Strong Landau damping:
\begin{equation}
\label{landau}
f(x,v,t=0)=\frac{1}{\sqrt{2\pi}}\left(1+\alpha\cos(kx)\right)\exp\left(-\frac{v^2}{2}\right),
\end{equation}
where $\alpha=0.5$, $k=0.5$.
\item  Two stream instability \Rmnum{1}:
\begin{equation}
\label{eq: two}
f(x,v,t=0)=\frac{2}{7\sqrt{2\pi}}(1+5v^2)\left(1+\alpha\left(\left(\cos(2kx)+\cos(3kx)\right)/1.2 + \cos(kx)\right)\right)\exp\left(-\frac{v^2}{2}\right),
\end{equation}
where $\alpha=0.01$, $k=0.5$.

\item Two stream instability \Rmnum{2}:
\begin{equation}
f(x,v,t=0)=\frac{1}{\sqrt{2\pi}}(1+\alpha\cos(kx))v^2\exp\left(-\frac{v^2}{2}\right),
\end{equation}
where $\alpha=0.05$, $k=0.5$.
\end{itemize}
The length of the domain in the $x{-}$
direction is $L=\frac{2\pi}{k}$ and the background ion distribution
function is fixed, uniform and chosen so that the total net charge
density for the system is zero.  To minimize the error from truncating the domain in the
$v$-direction, we let $v_{max} = 2\pi$. A fifth order SL finite difference WENO scheme is employed as a base scheme to achieve fifth order spatial accuracy.
Recall that the total mass, $L^p$ norms, entropy and total energy, which reads
\begin{align}\text{Mass}&=\int_v\int_xf(x,v,t)dxdv, \notag\\
\|f\|_{L^p}&=\left(\int_v\int_x|f(x,v,t)|^pdxdv\right)^\frac1p, \quad 1\leq p <\infty, \notag \\
\text{Entropy}& = \int_v\int_xf(x,v,t)\log(f(x,v,t))dxdv, \notag \\
\text{Energy}&=\frac12\left(\int_v\int_xf(x,v,t)v^2dxdv + \int_xE^2(x,t)dx\right), \notag
\end{align}
are conserved in the VP system. The proposed numerical schemes cannot exactly preserve these physical quantities except the total mass, yet it is of interest to track their time evolution to test the numerical schemes' performance.

Firstly, we validate the capacity of the proposed scheme to correct the low order splitting errors. We set the computational mesh as $N_x\times N_v=400\times400$. In the simulation, we fix the spatial mesh and compute the numerical solutions up to $T=0.1$ with different CFL numbers. A reference solution is computed with $\text{CFL} =0.01$ by using IDC3J3. In Table \ref{fig:landau}-\ref{fig:two2}, we report the $L^1$ error and the orders of accuracy when different orders of IDC methods coupled with first order splitting are used for the VP system with three sets of initial conditions. It is observed that the dimensional splitting error in time is significantly reduced when the IDC framework is applied, and $(K+1)^{th}$ order of accuracy is clearly achieved for IDC3J$K$ ($K\leq3$). In Table \ref{fig:idc_str}, the $L^1$ error and the orders of accuracy for IDC methods coupled with Strang splitting are reported, where $(2K+2)^{th}$ order of accuracy is observed for IDC-Strang3J$K$ ($K\leq1$). We remark that (1) In terms of stability, the Strang splitting is advantageous as it is known to be unconditionally stable, whereas the proposed scheme 
with IDC3J1 is only conditionally stable, see Table 4.1. When the CFL number is small enough, e.g., CFL = 0.6, 
the magnitudes of errors from IDC3J1 and the second order Strang splitting (IDC-Strang3J0 in Table 5.5) are comparable; 
however, the computational cost from the Strang splitting scheme is less than that of the IDC3J1 scheme; (2) the performance of IDC methods with first order splitting and Strang splitting are very similar when the same order accuracy and CFL numbers are considered. Hence, we only report the numerical results from IDC methods with first order splitting below for brevity.


\begin{table}[htb]
\begin{center}
\caption{Strong Landau damping. $T=0.1$. $L^1$ error and orders of accuracy.
\label{fig:landau}
}
\bigskip
\begin{tabular}{|c | c c|c c|c c|c c|}
\hline
\cline{1-9} &\multicolumn{2}{c|}{IDC3J0} &\multicolumn{2}{c|}{IDC3J1} &\multicolumn{2}{c|}{IDC3J2}  &\multicolumn{2}{c|}{IDC3J3}  \\
\cline{2-9} CFL& $L^1$ error&order&  $L^1$ error&order&  $L^1$ error&order&  $L^1$ error&order\\
\hline
0.6 &3.85E-06&	--       &9.71E-09	&--          &	2.22E-11&  --       &		1.83E-13&--	\\
\hline
0.5	&3.23E-06&	0.97&	6.79E-09&	1.96&	1.30E-11&	2.95&	8.90E-14&	3.95\\
\hline
0.4	&2.60E-06&	0.98&	4.37E-09&	1.97&	6.67E-12&	2.98&	3.66E-14&	3.98\\
\hline
0.3	&1.94E-06&	1.02&	2.44E-09&	2.03&	2.79E-12&	3.03&	1.15E-14&	4.03 \\
\hline
0.2	&1.30E-06&	0.99&	1.09E-09&	1.98&	8.34E-13&	2.98&	2.31E-15&	3.95 \\
\hline
\end{tabular}
\end{center}
\end{table}

\begin{table}[htb]
\begin{center}
\caption{Two stream instability \Rmnum{1}. $T=0.1$. $L^1$ error and orders of accuracy.
\label{fig:two1}}
\bigskip
\begin{tabular}{|c | c c|c c|c c|c c|}
\hline
\cline{1-9} &\multicolumn{2}{c|}{IDC3J0} &\multicolumn{2}{c|}{IDC3J1} &\multicolumn{2}{c|}{IDC3J2}  &\multicolumn{2}{c|}{IDC3J3}  \\
\cline{2-9} CFL& $L^1$ error&order&  $L^1$ error&order&  $L^1$ error&order&  $L^1$ error&order\\
\hline
0.6 &	4.46E-07&--		&2.36E-09&	--	&9.23E-12&	--	          &1.25E-13	&--\\
\hline
0.5	&3.73E-07	&0.98	&1.66E-09&	1.94&	5.43E-12&	2.91&	6.16E-14&	3.88\\
\hline
0.4	&2.99E-07	&1.00	&1.06E-09&	2.00&	2.78E-12&	3.00&	2.52E-14&	4.00\\
\hline
0.3	&2.24E-07	&1.00	&5.94E-10&	2.02&	1.17E-12&	3.02&	7.93E-15&	4.02\\
\hline
0.2	&1.50E-07	&0.99	&2.65E-10&	1.99&	3.48E-13&	2.98&	1.60E-15&	3.94\\
\hline

\hline
\end{tabular}
\end{center}
\end{table}

\begin{table}[htb]
\begin{center}
\caption{Two stream instability \Rmnum{2}. $T=0.1$. $L^1$ error and orders of accuracy.
\label{fig:two2}}
\bigskip
\begin{tabular}{|c | c c|c c|c c|c c|}
\hline
\cline{1-9} &\multicolumn{2}{c|}{IDC3J0} &\multicolumn{2}{c|}{IDC3J1} &\multicolumn{2}{c|}{IDC3J2}  &\multicolumn{2}{c|}{IDC3J3}  \\
\cline{2-9} CFL& $L^1$ error&order&  $L^1$ error&order&  $L^1$ error&order&  $L^1$ error&order\\
\hline
0.6&	6.72E-07&	--   &	2.33E-09&	--   &	3.61E-12&	--   &	1.79E-14&	--\\ \hline
0.5&	5.56E-07&	1.04&	1.61E-09&	2.04&	2.08E-12&	3.02&	8.60E-15&	4.01\\ \hline
0.4&	4.46E-07&	0.99&	1.03E-09&	2.00&	1.07E-12&	3.00&	3.53E-15&	3.99\\ \hline
0.3&	3.36E-07&	0.98&	5.84E-10&	1.98&	4.53E-13&	2.98&	1.13E-15&	3.95\\ \hline
0.2&	2.24E-07&	1.00&	2.60E-10&	2.00&	1.34E-13&	3.00&	2.65E-16&	3.58\\
\hline
\end{tabular}
\end{center}
\end{table}

\begin{table}[htb]
\begin{center}
\caption{IDC methods with Strang splitting for the VP systems. $T=0.1$. $L^1$ error and orders of accuracy.
\label{fig:idc_str}}
\bigskip
\begin{tabular}{|c | c c|c c|c c|c c|}
\hline
 &\multicolumn{4}{c|}{Strong Landau damping}  &\multicolumn{4}{c|}{Two stream instability \Rmnum{1}} \\
\cline{2-9} &\multicolumn{2}{c|}{IDC-Strang3J0} &\multicolumn{2}{c|}{IDC-Strang3J1} &\multicolumn{2}{c|}{IDC-Strang3J0}  &\multicolumn{2}{c|}{IDC-Strang3J1}  \\
\cline{2-9} CFL& $L^1$ error&order&  $L^1$ error&order&  $L^1$ error&order&  $L^1$ error&order\\
\hline
0.6	&5.27E-09&	--   & 	7.97E-14&	--	&1.07E-09	&	--&5.62E-14&--	\\ \hline
0.5	&3.68E-09&	1.96&	3.88E-14&	3.94&	7.55E-10&	1.91&	2.78E-14&	3.87\\ \hline
0.4	&2.37E-09&	1.98&	1.60E-14&	3.97&	4.89E-10&	1.94&	1.14E-14&	3.98\\ \hline
0.3	&1.32E-09&	2.03&	5.05E-15&	4.01&	2.82E-10&	1.91&	3.64E-15&	3.97\\ \hline
0.2	&5.93E-10&	1.98&	1.04E-15&	3.90&	1.37E-10&	1.78&	7.99E-16&	3.74\\

\hline
\end{tabular}
\end{center}
\end{table}

\begin{figure}[htb]
\begin{center}
\includegraphics[width=3in]{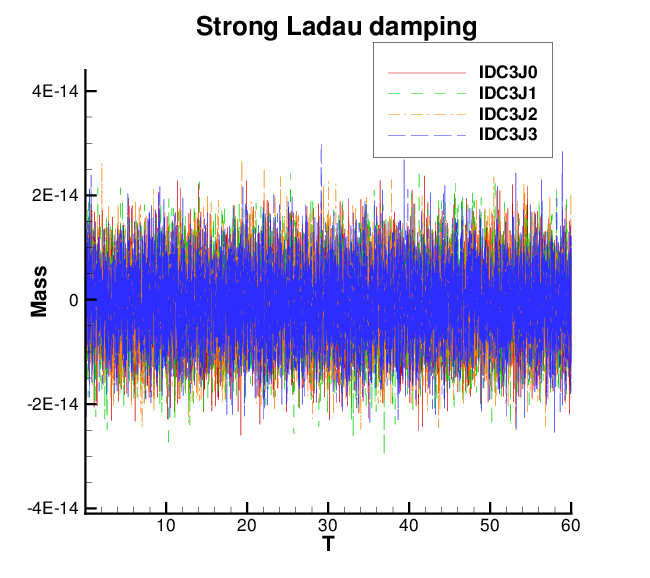}
\includegraphics[width=3in]{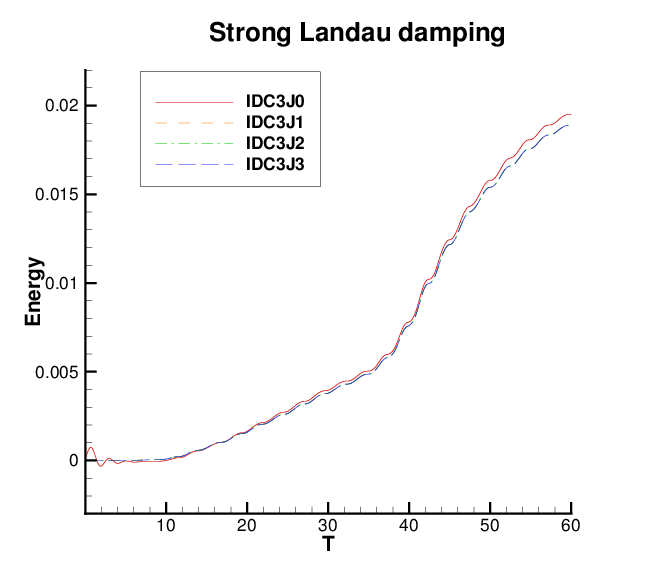}\\
\includegraphics[width=3in]{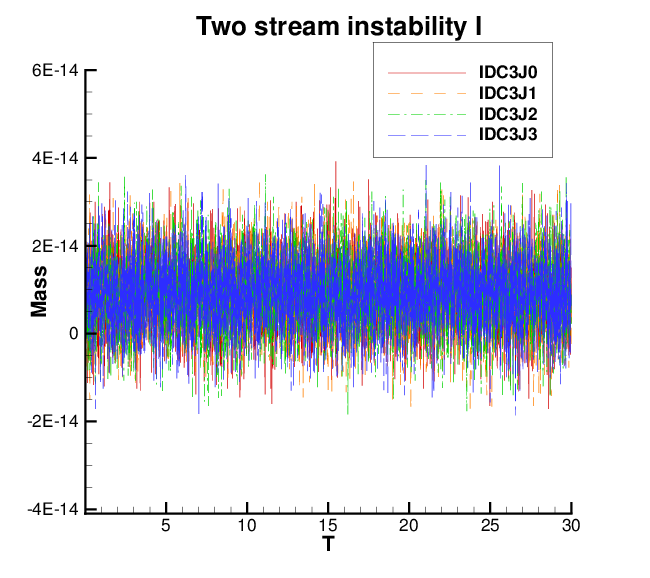}
\includegraphics[width=3in]{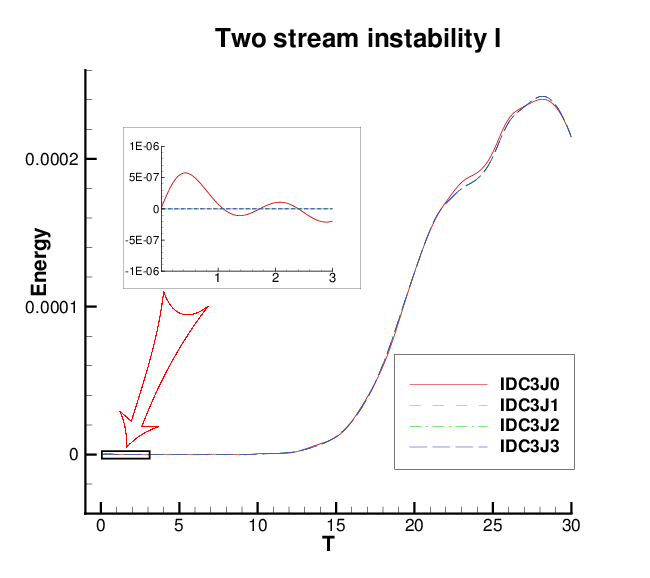}\\
\includegraphics[width=3in]{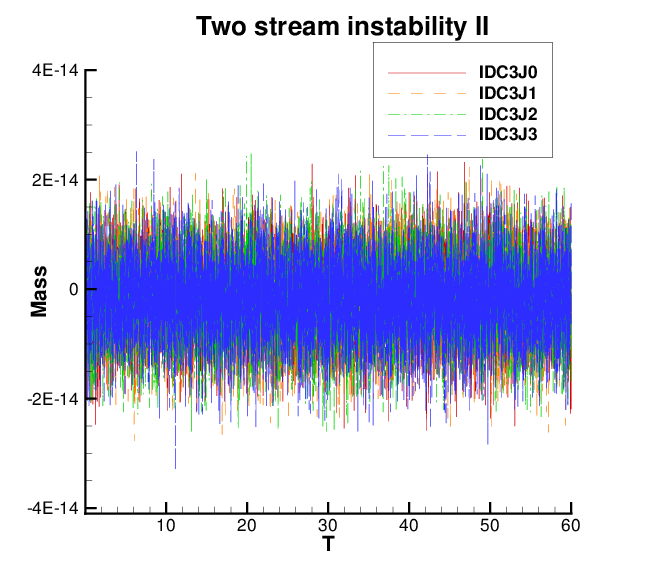}
\includegraphics[width=3in]{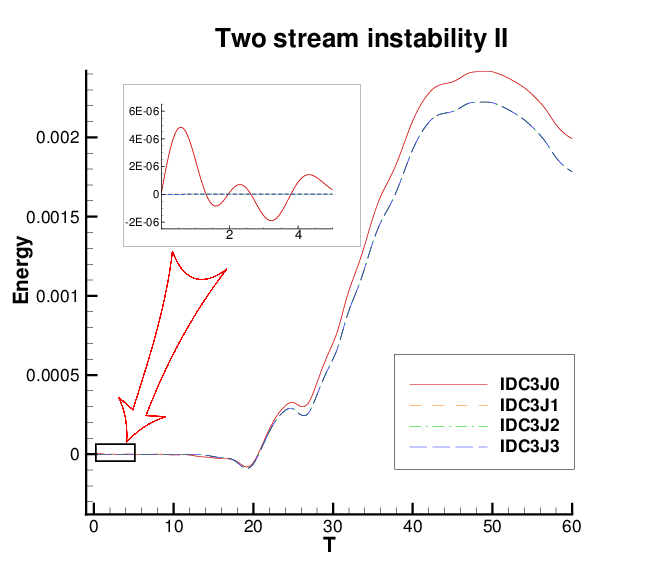}
\end{center}
\caption{
The time evolution of the relative deviation in total mass (left) and total energy (right). $N_x\times N_v=256\times256$. $\text{CFL} =0.6$.
}
\label{fig:time_evo}
\end{figure}

\begin{figure}[htb]
\begin{center}
\includegraphics[width=3in]{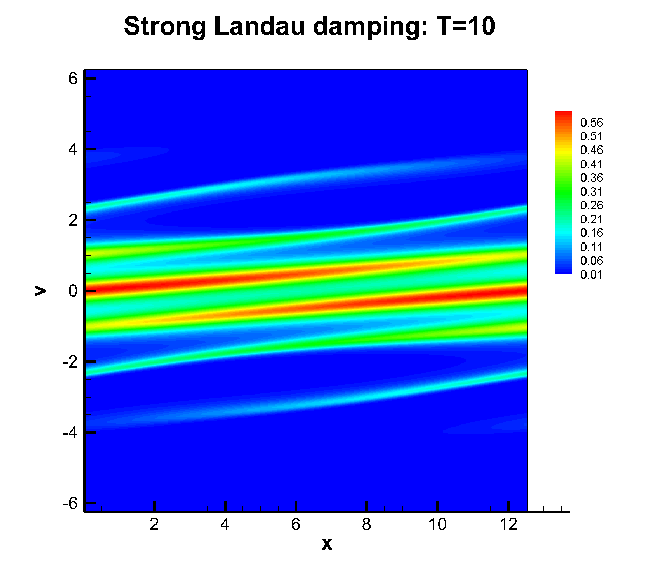}
\includegraphics[width=3in]{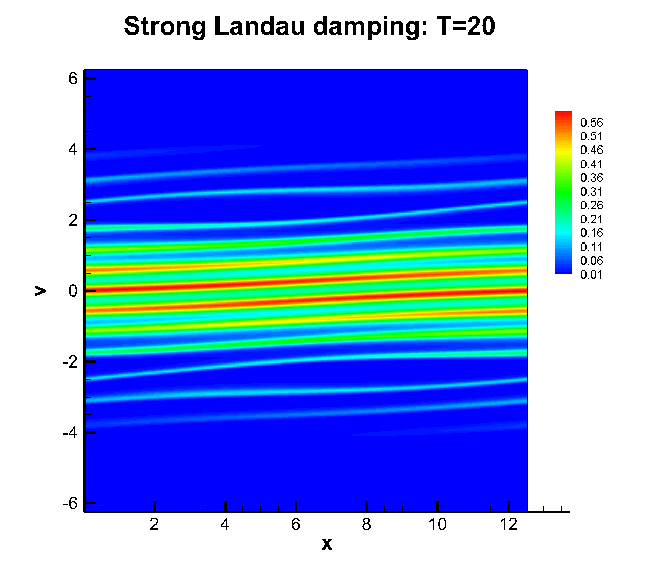}\\
\includegraphics[width=3in]{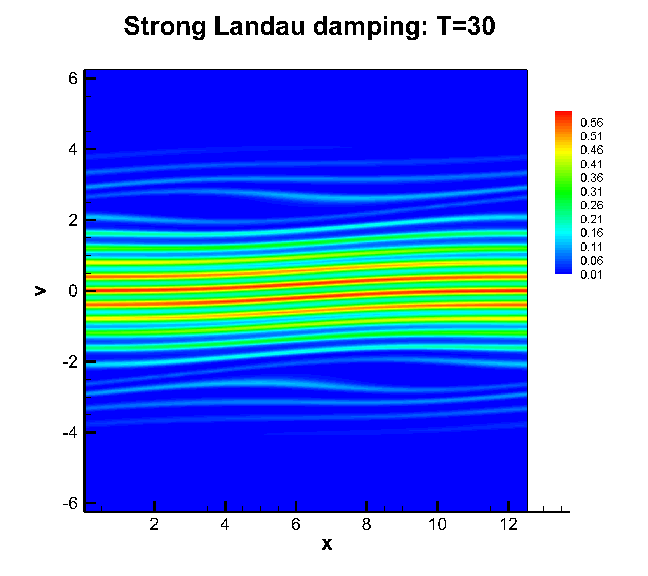}
\includegraphics[width=3in]{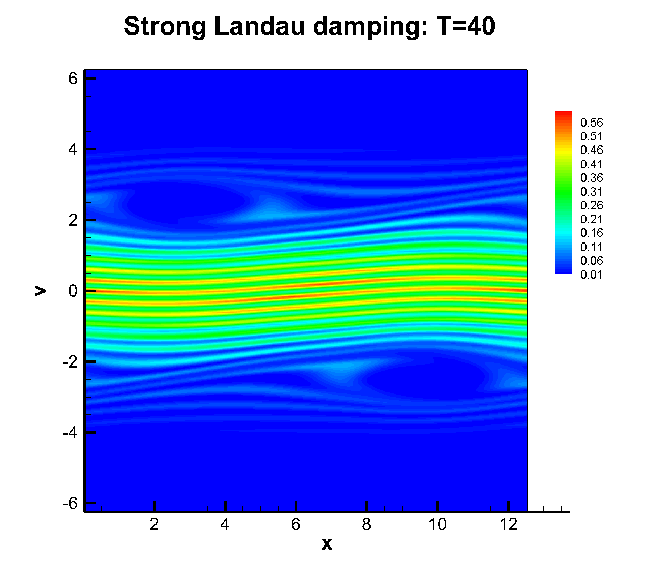}
\end{center}
\caption{
Contour plots of the numerical solutions for the strong Landau damping. $N_x\times N_v =256 \times 256$. $\text{CFL}  = 0.6$. IDC3J3.
}
\label{fig:lan_con}
\end{figure}

\begin{figure}[htb]
\begin{center}
\includegraphics[width=3in]{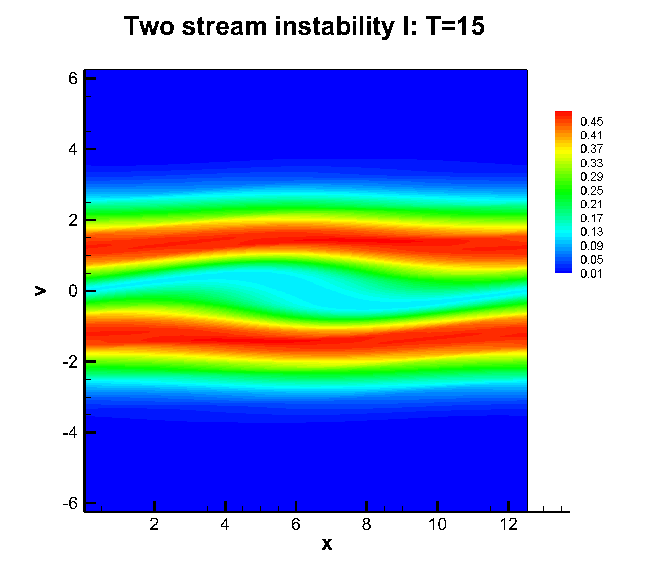}
\includegraphics[width=3in]{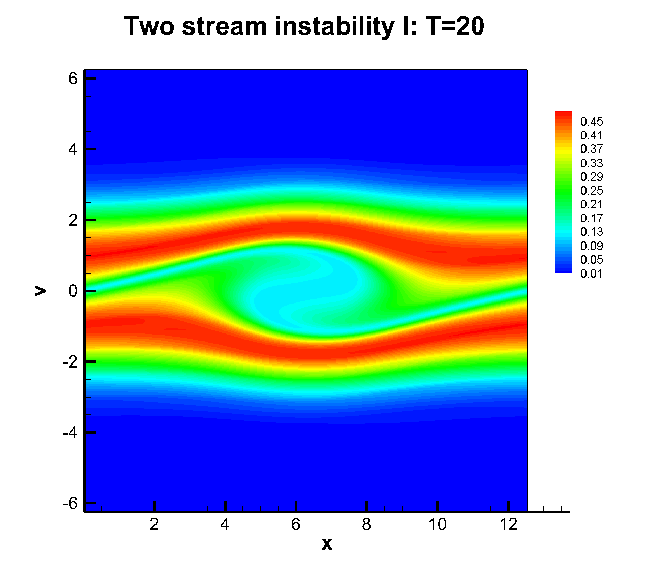}\\
\includegraphics[width=3in]{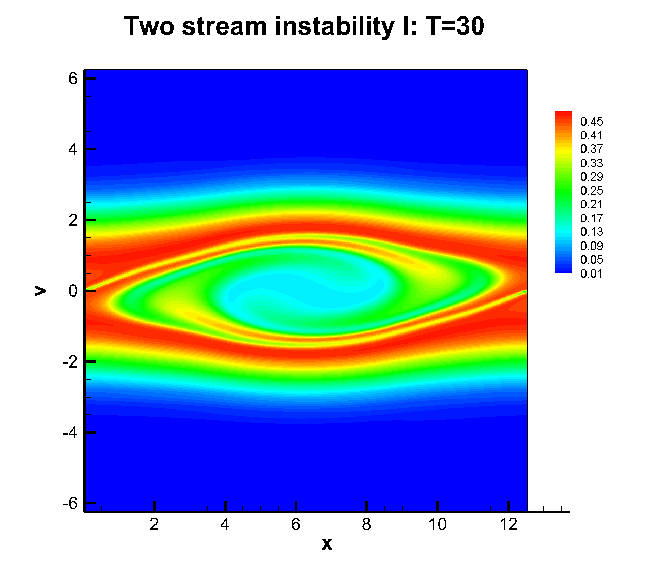}
\includegraphics[width=3in]{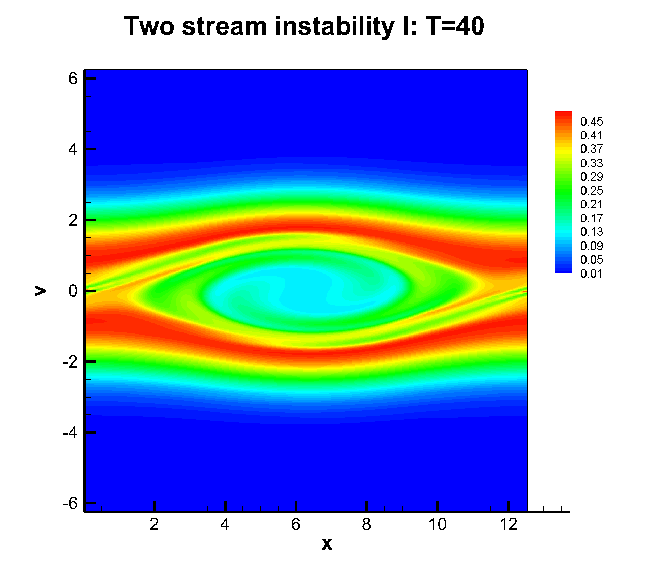}
\end{center}
\caption{
Contour plots of the numerical solutions for the two stream instability \Rmnum{1}. $N_x\times N_v =256 \times 256$. $\text{CFL}  = 0.6$. IDC3J3.
}
\label{fig:two_con1}
\end{figure}

\begin{figure}[htb]
\begin{center}
\includegraphics[width=3in]{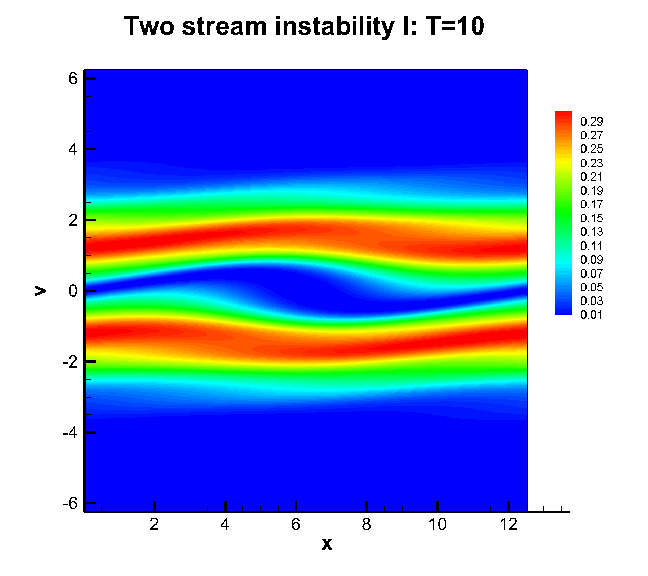}
\includegraphics[width=3in]{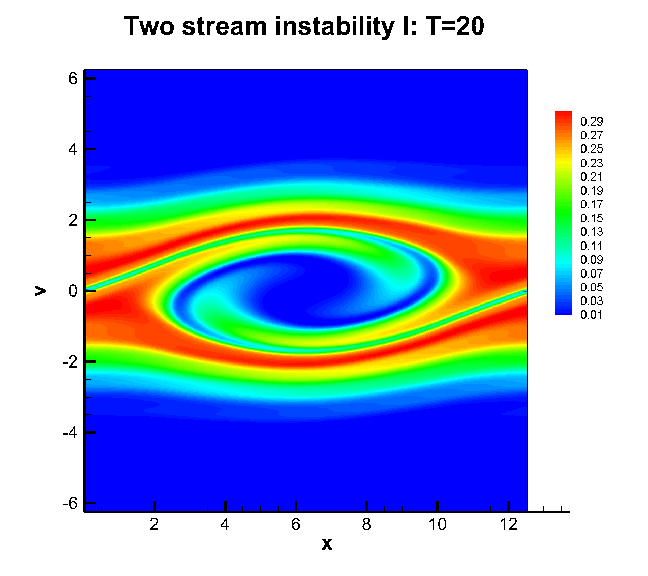}\\
\includegraphics[width=3in]{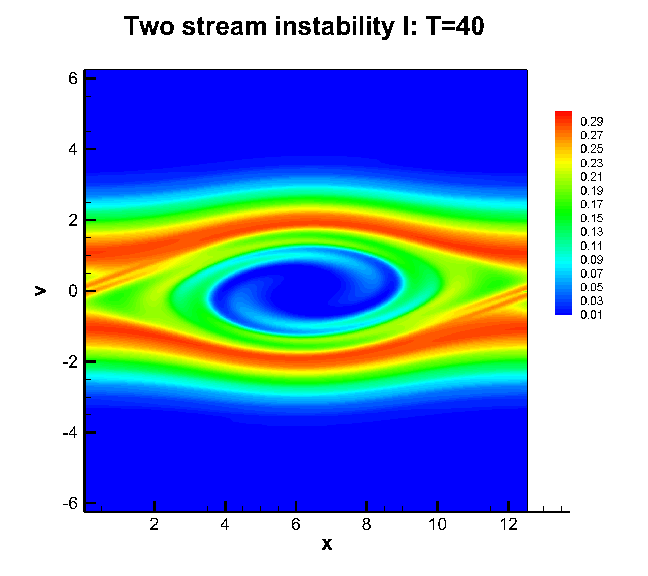}
\includegraphics[width=3in]{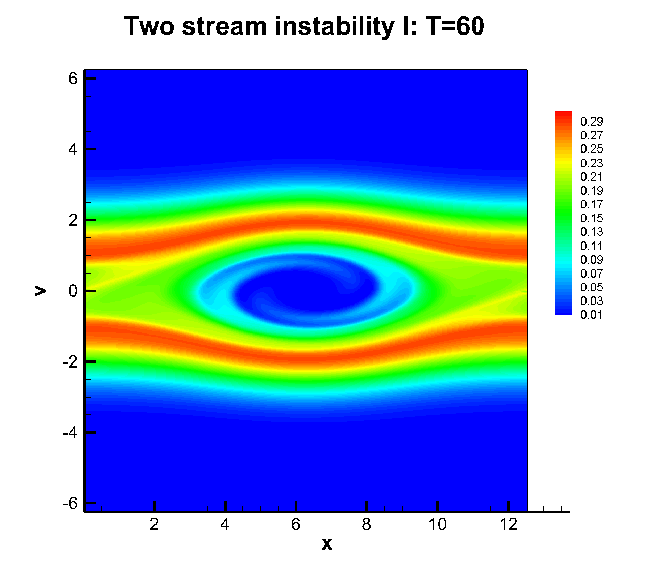}
\end{center}
\caption{
Contour plots of the numerical solutions for the two stream instability \Rmnum{2}. $N_x\times N_v =256 \times 256$. $\text{CFL}  = 0.6$. IDC3J3.
}
\label{fig:two_con}
\end{figure}

 Next we assess the ability of the schemes to preserve the above physical norms.
 The mesh is reset as $N_x\times N_v=256\times256$ and $\text{CFL} =0.6$. We track the evolution histories of these quantities and only show the results of total mass and energy here. In Figure \ref{fig:time_evo}, we report time evolution of the relative deviation in the total mass and the total energy for all three problems. It is observed that
the proposed scheme can preserve the total mass up to machine error as stated in Proposition \ref{prop: mass}. Moreover, the higher order schemes can better preserve the total energy than the first order scheme (IDC3J0). The performance of the second (IDC3J1), third (IDC3J2) and fourth (IDC3J3) order schemes are qualitatively similar. The reason may be that the error from the spatial discretization dominates the temporal error in these test cases. At last, we show the contour plots of the numerical solution for the strong Landau damping in Figure \ref{fig:lan_con}, two stream instability \Rmnum{1} and \Rmnum{2} in Figure \ref{fig:two_con1} and \ref{fig:two_con} respectively to demonstrate the performance of the proposed scheme. We only report the numerical results by the high order IDC3J3 method with $\text{CFL} =0.6$ for brevity. The results agree well with those presented in the literature by other methods \cite{heath2012discontinuous,James,Qiu_Christlieb,guo2012hybrid}.

\subsection{2-D Guiding center Vlasov equation}
\begin{figure}[htb]
\begin{center}
\includegraphics[width=3in]{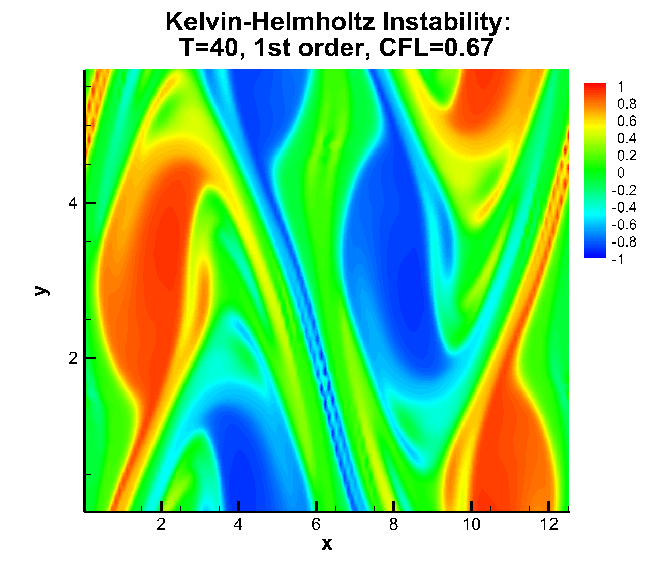}
\includegraphics[width=3in]{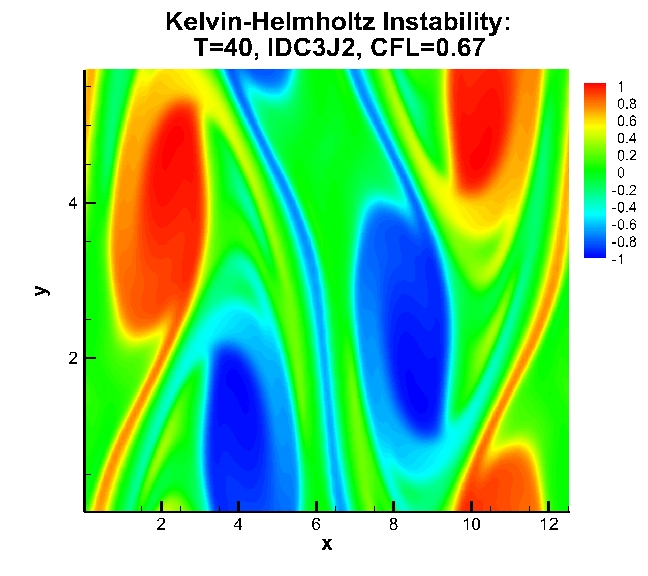}\\
\includegraphics[width=3in]{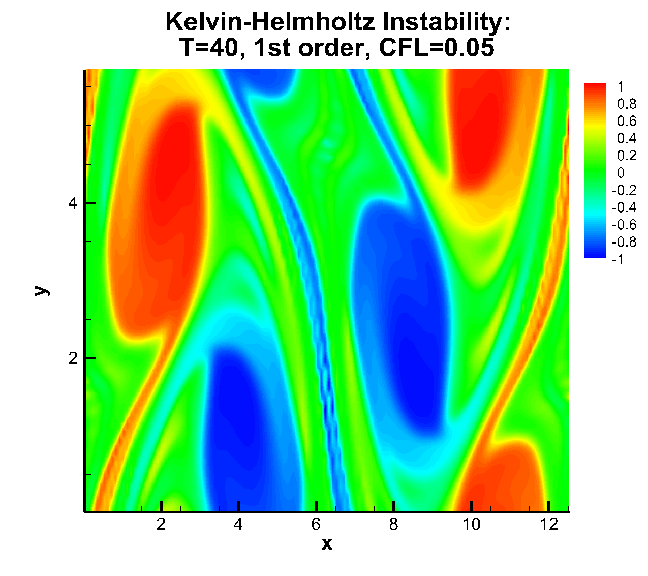}
\includegraphics[width=3in]{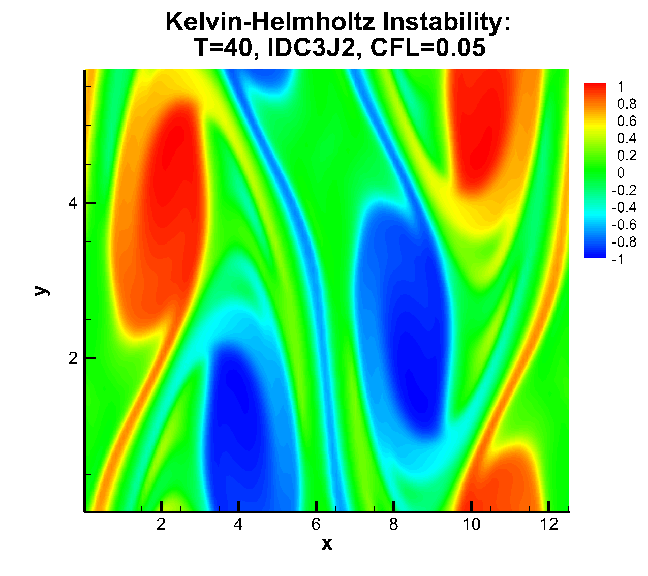}
\end{center}
\caption{
Contour plots of the numerical solutions for the Kelvin-Helmholtz instability. $N_x\times N_y =128 \times 128$. $\text{CFL}  = 0.67$ (top) and $\text{CFL} =0.05$ (bottom) at $T=40$. First order scheme (IDC3J0) (left); Third order scheme (IDC3J2) (right).
}
\label{fig:KH_con}
\end{figure}

\begin{figure}[htb]
\centerline{
\includegraphics[width=3in]{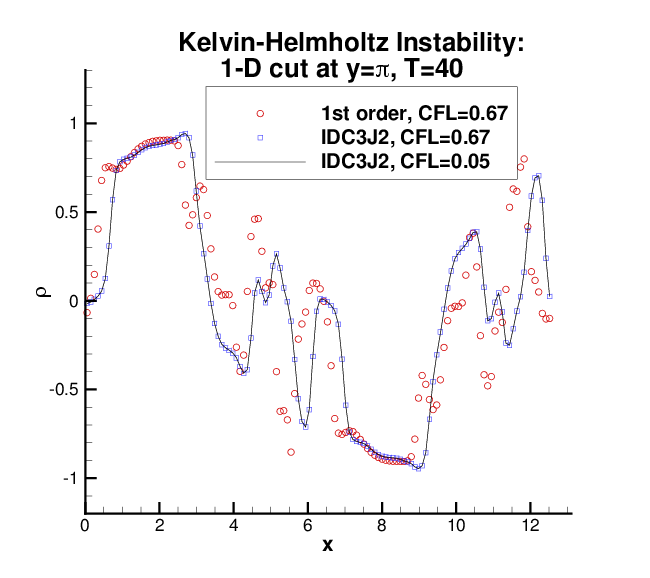}
\includegraphics[width=3in]{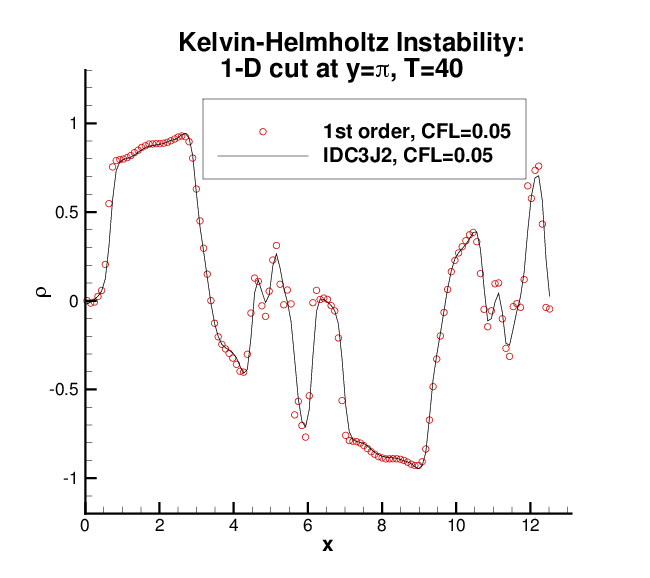}
}
\caption{1-D cuts of the numerical solutions at $y=\pi$ for the Kelvin-Helmholtz instability. $N_x\times N_y =128 \times 128$. $T=40$. $\text{CFL}  = 0.67$ (left); $\text{CFL}  = 0.05$ (right). The reference solution is computed by IDC3J2 with $\text{CFL} =0.05$.
}
\label{fig:KH_cut}
\end{figure}

\begin{figure}[htb]
\centerline{
\includegraphics[width=3in]{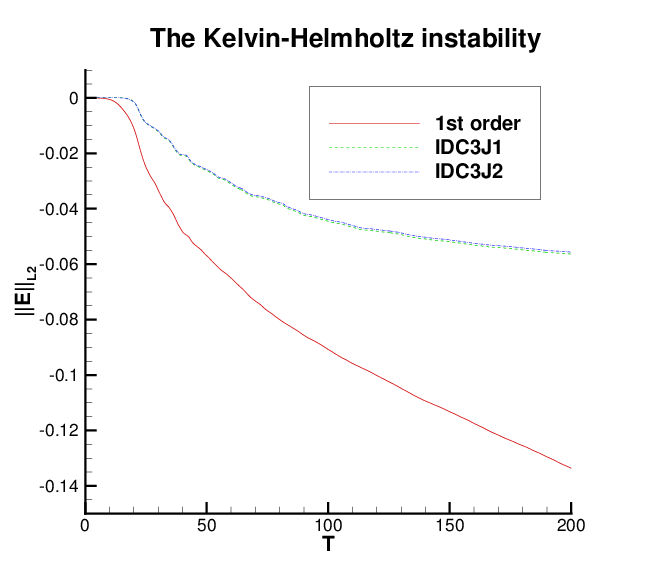}
\includegraphics[width=3in]{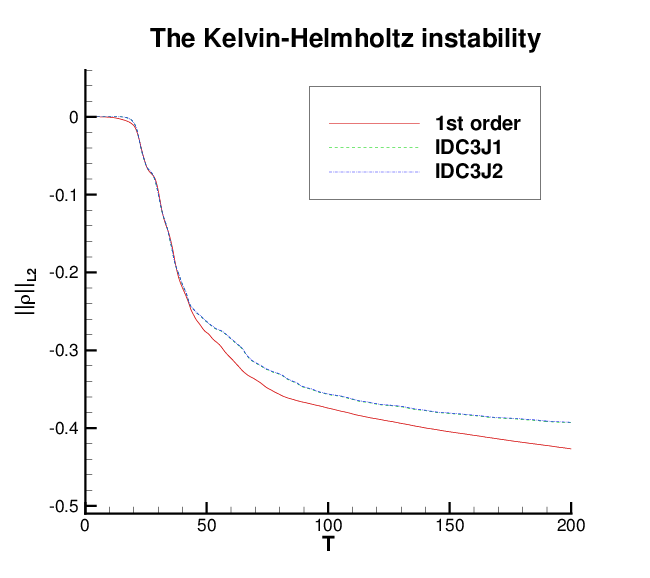}
}
\caption{The Kelvin-Helmholtz instability. The time evolutions of the energy $\|\mathbf{E}\|_{L^2}$ (left) and the enstrophy $\|\rho\|_{L^2}$ (right). $N_x\times N_y =128 \times 128$. $\text{CFL}  = 0.67$.
}
\label{fig:KH_evol}
\end{figure}

We consider the guiding center model \eqref{guiding} with
the initial condition \cite{crouseilles2010conservative}
$$\rho(x,y,0)=\sin(y)+ 0.015\cos(kx)$$
and periodic boundary conditions. We let $k=0.5$, thereby creating a Kelvin-Helmholtz instability. In the simulations, we set a mesh as $N_x\times N_y=128 \times 128$.
We use a third order SL finite difference WENO scheme as a base scheme to achieve third order spatial accuracy.
In Figure \ref{fig:KH_con} (top), we report the contour plots of numerical solutions at time $T=40$ for the first order splitting scheme (IDC3J0) and
the third order scheme (IDC3J2) with $\text{CFL} =0.67$. A noticeable difference is observed. In Figure \ref{fig:KH_cut} (left), 1-D cuts of the numerical solutions at $y=\pi$ are plotted. The reference solution
is computed by IDC3J2 with $\text{CFL} =0.05$. It can be observed that the numerical solution obtained by IDC3J2 with $\text{CFL} =0.67$ qualitatively matches the reference. However, a significant difference between
the numerical solution obtained by the first order scheme and the reference solution is observed. Then we reduce the CFL number to 0.05 for the first order splitting scheme. The 2-D
contour plot by the first order scheme approximately matches the reference by IDC3J2, see Figure \ref{fig:KH_con} (bottom). A more precise match is also observed when the 1-D cuts of the numerical solutions are compared. The presented numerical evidence shows better performance from higher order numerical schemes in time. Note that the continuous guiding center model \eqref{guiding} preserves the $L^2$ norms of $\rho$  (enstrophy) and the $L^2$ norms of $\mathbf{E}$ (energy), i.e.:
\[\frac{d}{dt}\|\rho(t)\|_{L^2}=\frac{d}{dt}\|\mathbf{E}(t)\|_{L^2}=0.\]
We track the relative deviations of these invariants numerically as a measurement of the quality of numerical schemes. In Figure \ref{fig:KH_evol}, the time evolutions of the enstrophy and the energy for the first order splitting scheme (IDC3J0), the second order scheme (IDC3J1) and the third order scheme (IDC3J2) are reported. It is observed that these quantities are better preserved by higher order schemes in time. Also note that little difference can be observed between the IDC3J1 and IDC3J2 solutions. In this case, the spatial error may dominate.

\subsection{2-D incompressible Euler equations}

We consider 2-D incompressible Euler equations in the vorticity stream-function formulation as follows
\beq
\label{eq: incomp_Euler}
\omega_t + \nabla \cdot ({\bf u} \omega) = 0, \quad x \in [0, 2\pi], \quad y\in[0, 2\pi].
\eeq
Here ${\bf u} = \nabla^\perp \Phi = (-\Phi_y, \Phi_x)$, where $\Phi$ is solved from the Poisson's equation.
We note that equation \eqref{eq: incomp_Euler} is in the same form as the 2-D guiding center Vlasov equation \eqref{guiding}.
We first test the accuracy of the scheme with the following initial conditions:
\begin{equation}
\omega(x,y,t=0) = -2\sin(x)\sin(y), \notag
\end{equation}
and periodic boundary conditions. Note that the exact solution is identical to the initial condition. We use this example to check the order of accuracy of our proposed scheme when the IDC framework is adopted to correct the splitting error. In the simulation, we fix the spatial mesh as $N_x\times N_y=300\times300$ and compute the numerical solutions with different CFLs. We evolve the solutions up to time $T=1$. In Table \ref{fig:euler}, we report the $L^1$ error and orders of accuracy for the first order splitting scheme (IDC3J0), the second order scheme (IDC3J1) and the third order scheme (IDC3J2). Expected orders of accuracy are observed. Note that third order convergence is not clearly observed for IDC3J2 when the CFL numbers are relatively small (CFL$\leq0.5$). In this case, the spatial error begins to dominate when the
numerical error is around $9.00E-9$. The $L^1$ error for IDC3J2 with $\text{CFL} =0.01$ is $9.31E-09$. We remark that there is a certain range of CFL numbers (which could be a very small interval), related to the spatial resolution and accuracy, where the temporal order of accuracy can be numerically observed. Above that range, the scheme is either numerically unstable or the order of convergence can not be observed yet; below that range, the spatial error could dominate, with which the temporal order of convergence can no longer be observed.

\begin{table} [htp]
\begin{center}
\caption{2-D incompressible Euler equation. A third order SL WENO scheme coupled in the IDC framework. $L^1$ norms of errors and orders of
accuracy. $N_x\times N_y=300\times300$. $T=1$.}
\smallskip
\begin{tabular}{|c|c c|c c|c c|
}\hline
\cline{2-7}
&\multicolumn{2}{c|}{IDC3J0}&\multicolumn{2}{c|}{IDC3J1}&\multicolumn{2}{c|}{IDC3J2}\\
\hline
 CFL & $L^1$ error & $L^1$ order  & $L^1$ error & $L^1$ order & $L^1$ error & $L^1$ order \\
\hline

0.67&	2.16E-03&	--&	8.89E-06	& --&	5.39E-08	&--  \\ \hline
0.62&	2.00E-03&	0.97&	7.63E-06&	1.97&	4.26E-08&	3.02 \\ \hline
0.57&	1.83E-03&	1.02&	6.43E-06&	2.03&	3.29E-08&	3.07 \\ \hline
0.52&	1.67E-03&	0.99&	5.36E-06&	1.99&	2.51E-08&	2.96 \\ \hline
0.47&	1.51E-03&	1.01&	4.37E-06&	2.01&	1.92E-08&	2.63 \\ \hline

\end{tabular}
\label{fig:euler}
\end{center}
\end{table}

\begin{figure}[htb]
\centerline{
\includegraphics[width=3in]{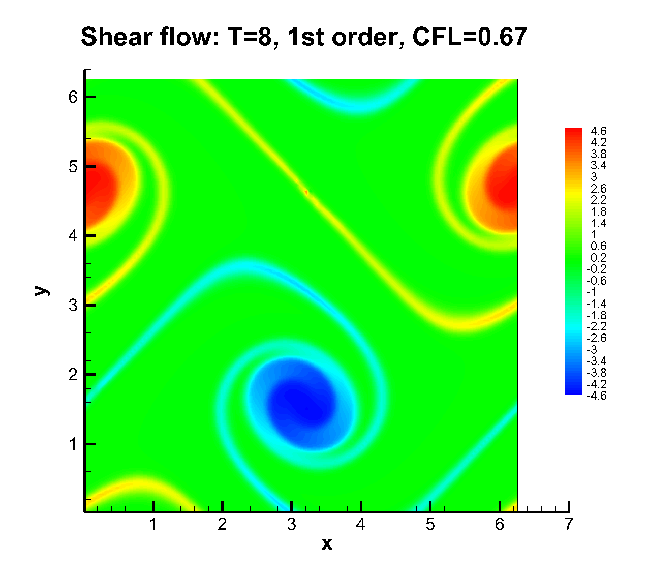}
\includegraphics[width=3in]{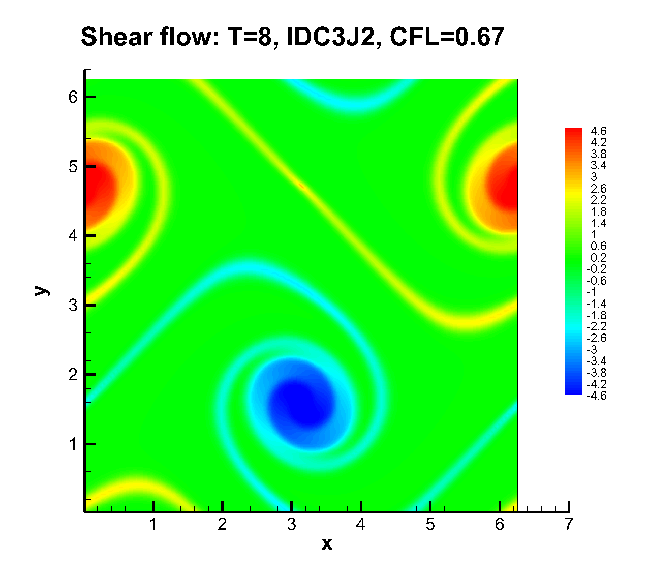}
}
\caption{
Contour plots of the numerical solutions for the shear flow test. $N_x\times N_y =128 \times 128$. $\text{CFL}  = 0.67$ at $T=8$. First order scheme (IDC3J0) (left); Third order scheme (IDC3J2) (right).
}
\label{fig:shear_con}
\end{figure}

\begin{figure}[htb]
\centerline{
\includegraphics[width=3in]{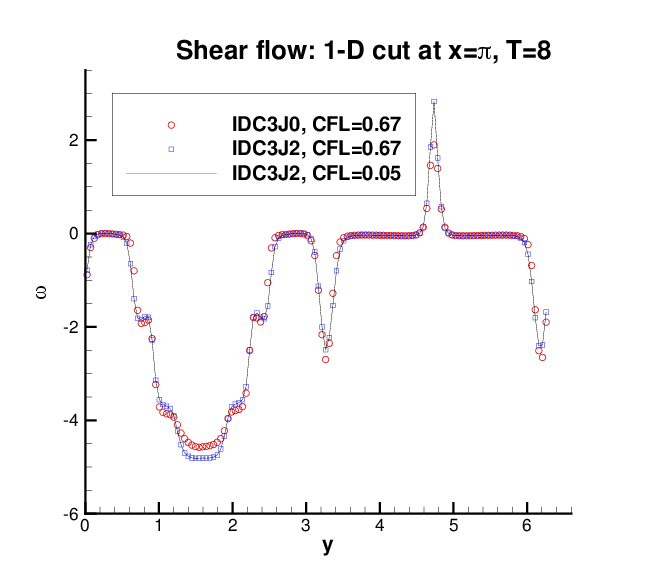}
\includegraphics[width=3in]{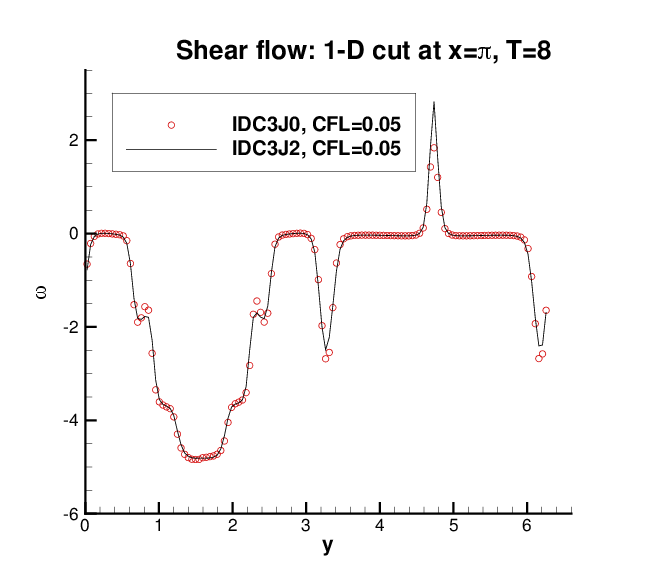}
}
\caption{1-D cuts of the numerical solutions at $x=\pi$ for the shear flow test. $N_x\times N_y =128 \times 128$. $T=8$. $\text{CFL}  = 0.67$ (left); $\text{CFL}  = 0.05$ (right). The reference solution is computed by IDC3J2 with $\text{CFL} =0.05$.
}
\label{fig:shear_cut}
\end{figure}

\begin{figure}[htb]
\begin{center}
\includegraphics[width=3in]{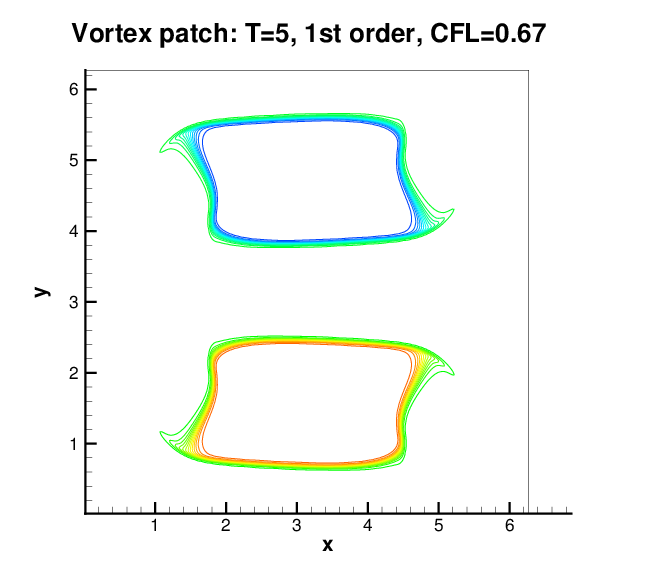}
\includegraphics[width=3in]{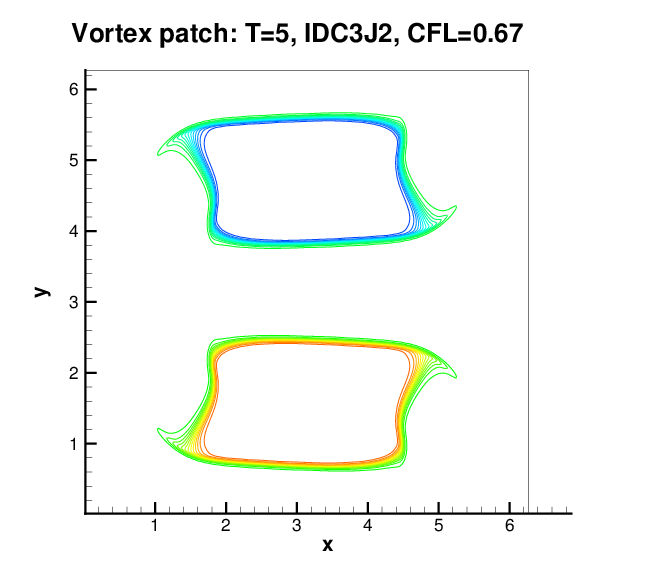}\\
\includegraphics[width=3in]{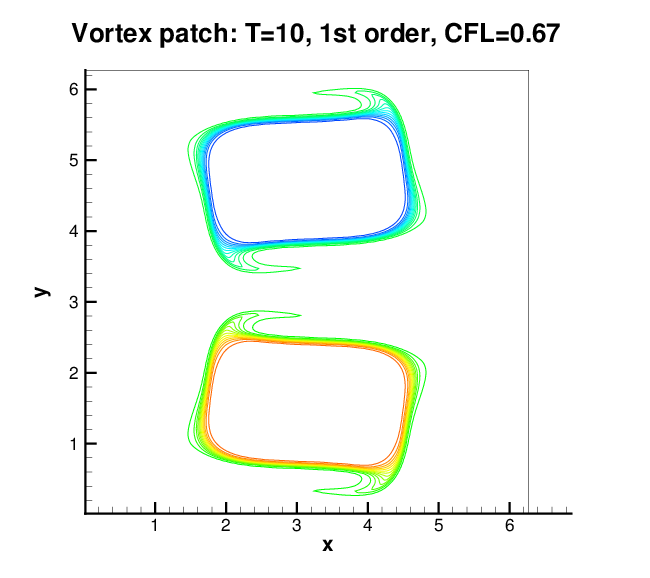}
\includegraphics[width=3in]{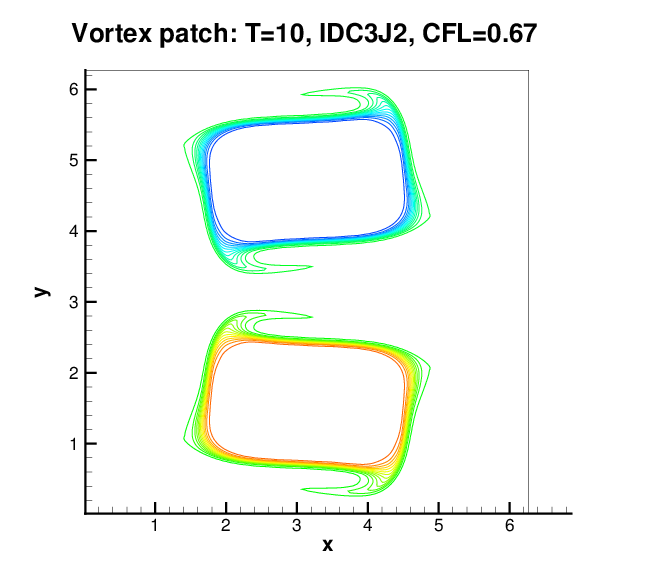}
\end{center}
\caption{
Contour plots of the numerical solutions for the vortex patch test. $N_x\times N_y =256 \times 256$. $\text{CFL}  = 0.67$ at $T=5$ (top) and $T=10$ (bottom). First order scheme (IDC3J0) (left); Third order scheme (IDC3J2) (right). 30 equally spaced contours from -1.1 to 1.1.
}
\label{fig:vor_con}
\end{figure}

We then consider two benchmark tests. One is shear flow with the initial condition given by:
\begin{equation}
\left\{\begin{array}{ll}
\displaystyle \omega(x, y, 0) = \delta \cos(x) - \frac{1}{\rho}\text{sech}^2((y-\pi/2)/\rho), & \text{if}\quad y\le \pi ,\\[5mm]
\displaystyle \omega(x, y, 0) = \delta \cos(x) + \frac{1}{\rho}\text{sech}^2((3\pi/2-y)/\rho), & \text{if}\quad y> \pi,
\end{array}\right.
\end{equation}
where $\delta = 0.05$ and $\rho = \frac{\pi}{15}$. The other is a vortex patch, and the initial condition is given by:
\begin{equation}
\left\{\begin{array}{ll}
\omega(x, y, 0) =-1, & \displaystyle\text{if}\quad (x, y) \in[\frac{\pi}{2},\frac{3\pi}{2}] \times [\frac{\pi}{4}, \frac{3\pi}{4}],\\[5mm]
\omega(x, y, 0) =1, &  \displaystyle\text{if}\quad (x, y)\in[\frac{\pi}{2}, \frac{3\pi}{2}]\times[\frac{5\pi}{4}, \frac{7\pi}{4}], \\[5mm]
\omega(x, y, 0) =0, & \text{ otherwise}.
\end{array}\right.
\end{equation}
In the simulations for the shear flow, we set the mesh as $N_x\times N_y = 128\times128$. A third order SL WENO scheme \cite{qiu_shu_sl} is used to obtain solutions in the IDC prediction step. In Figure \ref{fig:shear_con}, we report the contours of the numerical solution at time $T=8$ for the first order splitting scheme and IDC3J2. Little difference can be observed from the contour plot. To better see the difference, in Figure \ref{fig:shear_cut}, 1-D cuts at $x=\pi$ of the numerical solution are reported. We use the solution computed by IDC3J2 with $\text{CFL} =0.05$ as a reference. It is observed that the solution by IDC3J2 with $\text{CFL} =0.67$ matches the reference solution very well, whereas there is a noticeable difference between the solution by the first order scheme with $\text{CFL} =0.67$ and the reference. Then we reduce the CFL to 0.05 for the first order scheme and the corresponding temporal error is reduced. A more precise match is observed. The comparison shows the better performance of coupling the high order IDC methods to correct low order splitting errors. For the vortex patch test, we set the mesh size as $N_x\times N_y = 256\times256$. Figure \ref{fig:vor_con} gives the numerical solutions at time $T=5$ (top) and $T=10$ (bottom) for the first order scheme (IDC3J0) (left) and IDC3J2 (right). The solution structure is observed to be slightly better resolved when the high order IDC framework is used.

\section{Conclusion}
\label{conclusion}

In this paper, we proposed numerical methods that are high order in space and in time for Vlasov simulations. The proposed methods couple the SL finite difference WENO scheme with the IDC framework: we adopt the dimensional splitting SL WENO scheme as a base scheme to get a predicted solution in IDC, and the low order dimensional splitting error is iteratively reduced by solving the error equations again in a  dimensional splitting fashion in the correction steps in IDC.
We extended the scheme to solve the guiding center Vlasov equation and the two dimensional incompressible flow in vorticity stream-function formulation. A collection of numerical experiments demonstrate great performance of the proposed high order scheme in resolving solutions' structures, even in the long term. Unfortunately, the IDC framework renders some CFL time step restriction, despite the SL evolution mechanism in the prediction and correction steps of IDC. We quantify such CFL restrictions via Fourier analysis for several high order methods that we used in simulations. Another potential difficulty is the implementation of boundary conditions in a high order manner for problems with non-periodic boundary conditions. These are subject to future investigations.

\bigskip
\noindent
{\bf Acknowledgement.} The authors thank Dr. Matt Causley  for many helpful comments and suggestions.

\bibliographystyle{siam}
\bibliography{refer}

\end{document}